\documentclass[preprint,12pt]{elsarticle}
\usepackage{}

\usepackage{amsmath}
\usepackage{cases}
\usepackage{color}
\usepackage{hyperref}
\usepackage{amssymb}
\usepackage{graphicx}
\usepackage{graphics}
\usepackage{subfigure}
\usepackage{appendix}
\usepackage{textcomp}
\usepackage{bm}
\usepackage{url}

\makeatletter
  \newcommand\figcaption{\def\@captype{figure}\caption}
  \newcommand\tabcaption{\def\@captype{table}\caption}
\makeatother

\biboptions{numbers,sort&compress}

\begin{document}

\begin{center}

{\Large {Approximate Solutions of Von K{\'a}rm{\'a}n Plate}  {under Uniform Pressure} 
--- Equations in Differential Form}

\vspace{0.3cm}

 Xiaoxu Zhong $^3$, Shijun Liao $^{1,2,3}$  \footnote{Corresponding author.  Email address: sjliao@sjtu.edu.cn}

\vspace{0.3cm}

$^1$ State Key Laboratory of Ocean Engineering, Shanghai 200240, China\\
$^2$ Collaborative Innovative Center for Advanced Ship and Deep-Sea Exploration, Shanghai 200240, China\\
$^3$ School of Naval Architecture, Ocean and Civil Engineering\\
Shanghai Jiao Tong University, Shanghai 200240, China

 \end{center}

\hspace{-0.6cm}{\bf Abstract}  {\em The large deflection of a circular thin plate under uniform external pressure is a classic problem in solid mechanics, dated back to Von K{\'a}rm{\'a}n \cite{Karman}.  {This problem is reconsidered in this paper using an analytic approximation method, namely the homotopy analysis method (HAM).}  Convergent series solutions are obtained for four types of boundary conditions with rather high nonlinearity, even in the case of $w(0)/h>20$, where $w(0)/h$ denotes  the ratio of central deflection to plate thickness.  Especially,  we prove that the  previous perturbation methods for an arbitrary  perturbation quantity (including the Vincent's  \cite{Vincent} and Chien's \cite{Qian} methods) and the modified iteration method \cite{YehLiu}  are only the special cases of the HAM.  However, the HAM  works well even when the perturbation methods become invalid.   All of these demonstrate the validity and potential of the HAM for the Von K{\'a}rm{\'a}n's  plate equations, and  show the superiority of the HAM over perturbation  methods for highly nonlinear problems. }

\vspace{0.3cm}

\hspace{-0.6cm}{\bf Key Words} circular plate, \sep uniform external pressure, \sep homotopy analysis method (HAM)

\section{Introduction}

The  large deflection of a circular thin plate under uniform external pressure, which may date back to Love \cite{Love} and  shell \cite{VonTsien},  plays an important role in many engineering fields, such as mechanical and marine engineering, the precision instrument manufacture, and so on.  In 1910,  Von K{\'a}rm{\'a}n \cite{Karman} derived the so-called Von K{\'a}rm{\'a}n's plate equations for large displacements.     Its differential form \cite{Karman,Zheng}  reads
\begin{eqnarray}
   {\cal N}_1[\varphi,S,y] &=& y^{2}\frac{d^{2}\varphi(y)}{dy^{2}}-\varphi(y)S(y)-Qy^{2}=0, \label{geq1:original:U}\\
    {\cal N}_2[\varphi,S,y] &=& y^{2}\frac{d^{2}S(y)}{dy^{2}}+\frac{1}{2}\varphi^{2}(y)=0,
    \label{geq2:original:U}
\end{eqnarray}
subject to the boundary conditions
\begin{equation}
    \varphi(0)=S(0)=0,  \label{ic1:original:U}
\end{equation}
\begin{equation}
    \varphi(1)=\frac{\lambda}{\lambda-1}\cdot\frac{d\varphi(y)}{dy}\bigg{|}_{y=1},\;\;\;\;S(1)=\frac{\mu}{\mu-1}\cdot\frac{dS(y)}{dy}\bigg{|}_{y=1},  \label{ic2:original:U}
\end{equation}
under the definitions
\begin{equation}
    y=\frac{r^{2}}{R_{a}^{2}},\;\;\;\;W(y)=\sqrt{3(1-\nu^{2})}\frac{w(y)}{h},\;\;\;\;\varphi(y)=y\frac{dW(y)}{dy},  \label{dimensionless1}
\end{equation}
\begin{equation}
    S(y)=3(1-\nu^{2})\frac{R_{a}^{2}N_{r}}{Eh^{3}}y,\;\;\;\;Q=\frac{3(1-\nu^{2})\sqrt{3(1-\nu^{2})}R_{a}^{4}}{4Eh^{4}}p, \label{dimensionless2}
\end{equation}
where $r$ is the radial coordinate whose origin locates at the center of the plate, $w(y)$ and $N_{r}$ describe the deflection and the radial membrane force of the plate, the constants $E$, $\nu$, $R_{a}$, $h$ are elastic modulus, the Poisson's ratio, radius and thickness of the plate, respectively, $p$ represents the external uniform load, $\lambda$ and $\mu$ are parameters related to the boundary conditions at $y=1$.  From Eq.~(\ref{dimensionless1}),  we have the dimensionless central deflection
\begin{equation}
    W(y)=-\int_{y}^{1}\frac{1}{\varepsilon}\varphi(\varepsilon)d\varepsilon.   \label{deflection}
\end{equation}

\begin{figure}[!htb]
    \centering
    \subfigure[Clamped]{
    \label{clamped:figure}
    \includegraphics[width=2.3in]{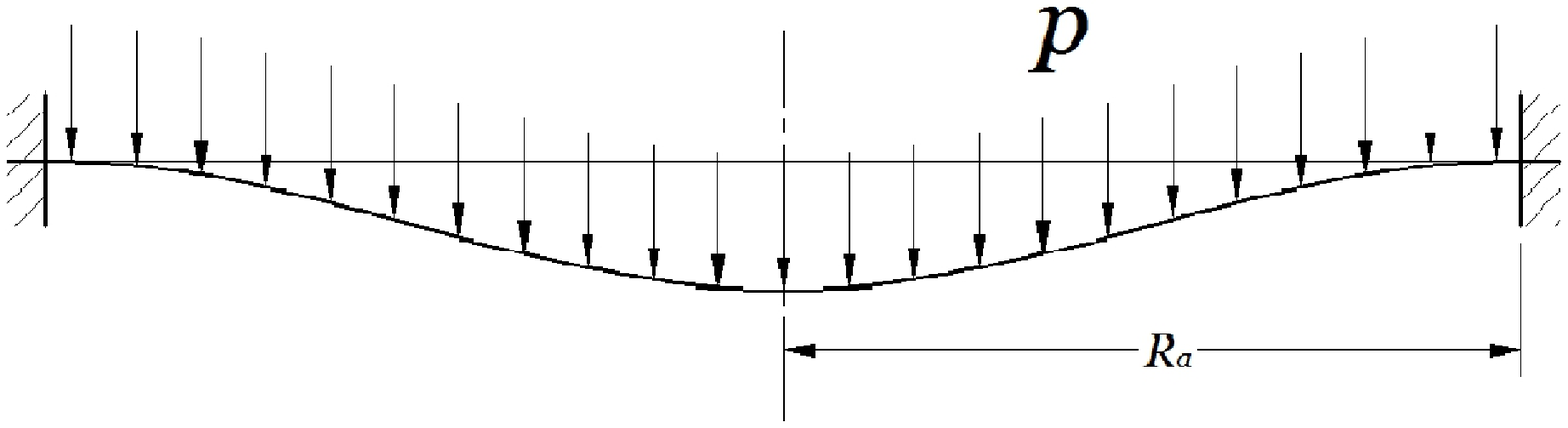}}
    \hfill
    \subfigure[Moveable clamped]{
    \label{moveable:clamped:figure}
    \includegraphics[width=2.3in]{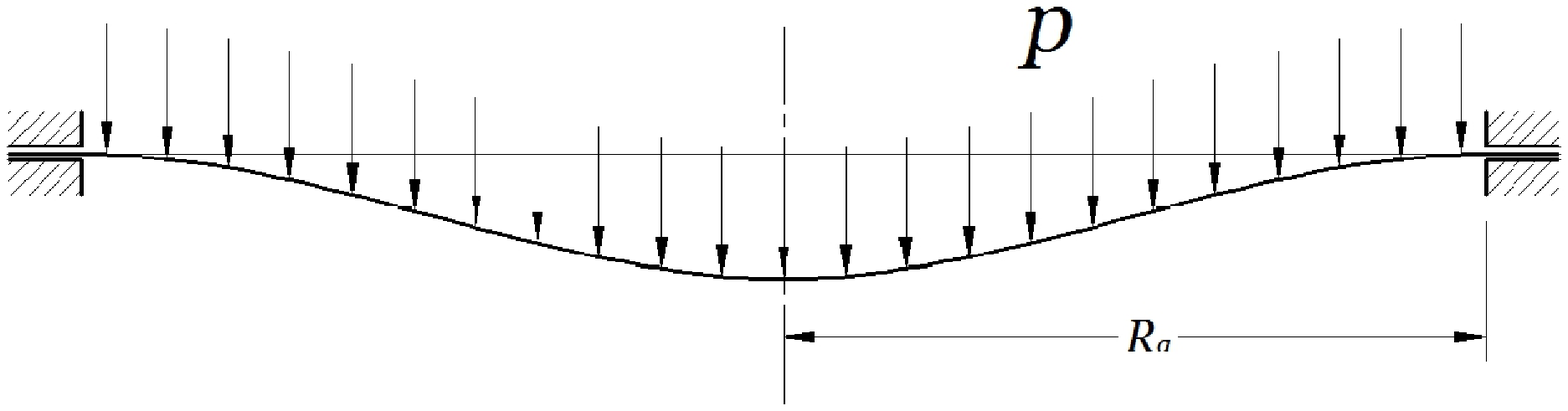}}
   \vfill
   \subfigure[Simple support]{
    \label{simple:supported:figure}
    \includegraphics[width=2.3in]{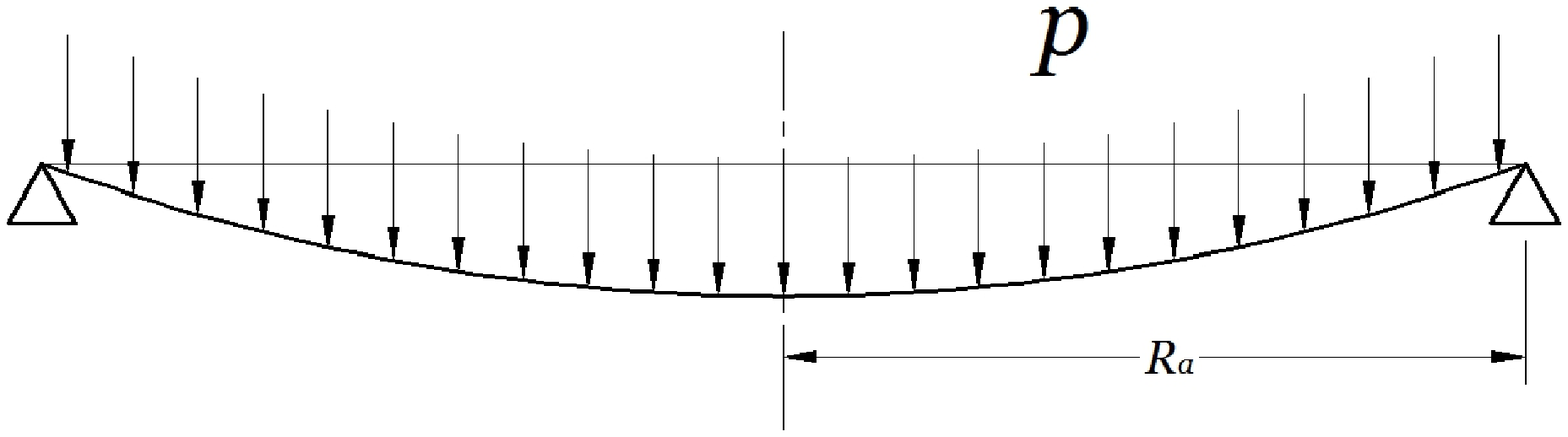}}
    \hfill
    \subfigure[Simple hinged support]{
    \label{simple:hinged:figure}
    \includegraphics[width=2.3in]{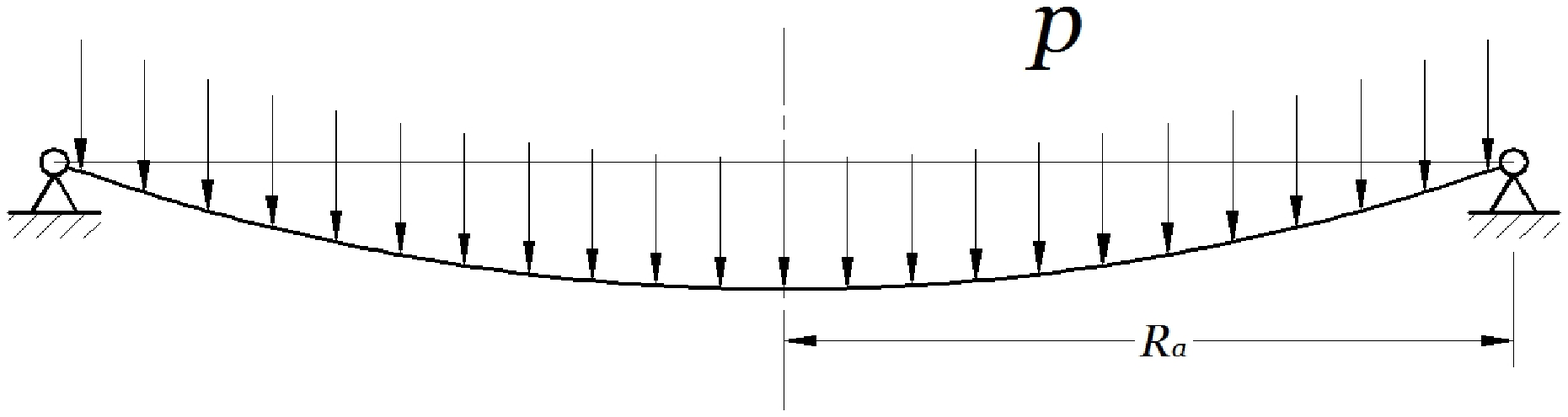}}
    \caption{{Four types of boundary conditions}}
    \label{four:types:figure}
\end{figure}

 {As shown in Fig.~\ref{four:types:figure}}, four types of boundary conditions are considered:
 \begin{enumerate}
\item[(a)] Clamped: $\lambda=0$ and $\mu=2/(1-\nu)$;
\item[(b)] Moveable clamped: $\lambda=0$ and $\mu=0$;
\item[(c)] Simple support: $\lambda=2/(1+\nu)$ and $\mu=0$;
\item[(d)] Simple hinged support: $\lambda=2/(1+\nu)$ and $\mu=2/(1-\nu)$;
\end{enumerate}

 Over the past century, lots of analytic/numerical  methods are proposed \cite{Way, Chien2,DaDeppo, Simmonds, Zhou, Jian, ZhengLee, Mahmoud} for the Von K{\'a}rm{\'a}n's plate equations.  Vincent \cite{Vincent} proposed a perturbation approach by using the load  $Q$  as a small physical parameter,  say, 
  \begin{equation}
 \begin{split}
    \varphi(y)=\varphi^{(V)}(y)=\sum_{i=1}^{+\infty}\varphi_{i}^{(V)}(y)Q^{2i-1},\\
    S(y)=S^{(V)}(y)=\sum_{i=1}^{+\infty}S_{i}^{(V)}(y)Q^{2i},\;\;\;
    \label{Vincent:Expand:S}
 \end{split}
\end{equation}
where $\varphi_{i}^{(V)}(y)$ and $S_{i}^{(V)}(y)$ are governed by 
\begin{equation}
{Q:} \left\{
\begin{split}
{y^2\frac{d^{2}\varphi_{1}^{(V)}(y)}{dy^{2}} = }{~y^{2}} {, \;\;\;\;y\in(0,1), }~~~~~~~~ ~~~~~~~~~~~~~~~~~~~~~~~ ~\label{Vincent:Q1}\\
 {\varphi_{1}^{(V)}(0)=0,\;\;\;\;\;\varphi_{1}^{(V)}(1)=\frac{\lambda}{\lambda-1}\frac{d\varphi_{1}^{(V)}(y)}{dy}\bigg{|}_{y=1};\;\;\;\;~~~~~~~~~~}
\end{split}
\right.
\end{equation}
\begin{equation}
 {Q^{2}:} \left\{
\begin{split}
{y^2\frac{d^{2}S_{1}^{(V)}(y)}{dy^{2}}=-\frac{1}{2}\left(\varphi_{1}^{(V)}(y)\right)^2, } ~~~~~~~~~~~~~~~~ ~~~~~~~~~~~~~\\
{S_{1}^{(V)}(0)=0,\;\;\;\;\;S_{1}^{(V)}(1)=\frac{\mu}{\mu-1}\frac{dS_{1}^{(V)}(y)}{dy}\bigg{|}_{y=1};~~~~~~~~~~}
\end{split}
\right.
\end{equation}
{$~~~~~\cdots \cdots $}
\begin{equation}
 {Q^{2i+1}:} \left\{
\begin{split}
 {y^2\frac{d^{2}\varphi_{i+1}^{(V)}(y)}{dy^{2}}=\sum_{j=1}^{i}\varphi_{j}^{(V)}(y)S_{i-j+1}^{(V)}(y), } ~~~~~~~~~~~~~~~~~~~~~~~~~\\
 {\varphi_{i+1}^{(V)}(0)=0,\;\;\;\;\varphi_{i+1}^{(V)}(1)=\frac{\lambda}{\lambda-1}\frac{d\varphi_{i+1}^{(V)}(y)}{dy}\bigg{|}_{y=1};~~~~~~~~~~~~~~}
\end{split}
\right.
\end{equation}
\begin{equation}
 {Q^{2i+2}:} \left\{
\begin{split}
 {y^2\frac{d^{2}S_{i+1}^{(V)}(y)}{dy^{2}}=-\frac{1}{2}\sum_{j=1}^{i+1}\varphi_{j}^{(V)}(y)\varphi_{i-j+2}^{(V)}(y), ~~~~~~~~~~~~~~~~~~~~}\; \\
 {S_{i+1}^{(V)}(0)=0,\;\;\;\;\;S_{i+1}^{(V)}(1)=\frac{\mu}{\mu-1}\frac{dS_{i+1}^{(V)}(y)}{dy}\bigg{|}_{y=1}.~~~~~~~~~~~~~~}
\end{split}
\right. \label{Vincent:Q:universe}
\end{equation}
 { In the case of circular plate with clamped boundary ($\lambda=0$, $\mu=20/7$), the Vincent's perturbation solutions are:}
\begin{eqnarray}
 {\varphi_{1}^{(V)}(y)} &=&  {-\frac{y}{2}+\frac{y^2}{2}, } \label{Vincent:phi1}\\
 {S_{1}^{(V)}(y)} &=&  {\frac{41y}{672}-\frac{y^2}{16}+\frac{y^3}{24}-\frac{y^4}{96},}\\
 {\varphi_{2}^{(V)}(y)} &=&  {\frac{659y}{80640}-\frac{41y^2}{2688}+\frac{83y^3}{8064}-\frac{5y^4}{1152}+\frac{y^5}{768}-\frac{y^6}{5760},}\\
 {S_{2}^{(V)}(y)} &=&  {-\frac{2357y}{1505280}+\frac{659y^2}{322560}-\frac{1889y^3}{967680}+\frac{103y^4}{96768}-\frac{59y^5}{161280} } \nonumber\\
&\;\;& {+\frac{13y^6}{138240}-\frac{17y^7}{967680}+\frac{y^8}{645120},} \label{Vincent:S2}\\
\cdots \cdots&\;\;& \nonumber
\end{eqnarray}

Unfortunately, Vincent's perturbation method  is valid only for rather small ratio of central deflection to plate thickness $w(0)/h<0.52$ for a circular plate with clamped boundary  \cite{Vincent}.     Thereafter, extensive researches were done to find a better perturbation quantity.   For instance,  the central deflection and the average angular deflection were used by Chien \cite{Qian}, Chien and Yeh \cite{QianYeh}  and Hu \cite{Hu} as perturbation quantity.  It is found \cite{Chen} that the central deflection is the best, which can give perturbation results convergent within $w(0)/h<2.44$ for clamped boundary  \cite{Qian}.   Expanding $\varphi(y)$, $S(y)$ and $Q$ into power series of the central deflection $W(0)$, we have 
\begin{equation}
\left\{
\begin{split}
    \varphi(y)=\varphi^{(C)}(y)=\sum_{m=1}^{+\infty}\varphi_{m}^{(C)}(y)W^{2m-1}(0),\\
S(y)=S^{(C)}(y)=\sum_{m=1}^{+\infty}S_{m}^{(C)}(y)W^{2m}(0),~~\\
Q=Q^{(C)}=\sum_{m=1}^{+\infty}Q_{m}^{(C)}W^{2m-1}(0).~~~~~~~
\label{Chien:Expand:Q}
\end{split} \right.
\end{equation}
 {The procedures of Chien's perturbation method \cite{Qian,Zheng} are:}
\begin{equation}
 {W(0):} \left\{
\begin{split}
 {y^2\frac{d^{2}\varphi_{1}^{(C)}(y)}{dy^{2}}=Q_{1}^{(C)}y^{2}, \;\;\;\;y\in(0,1), ~~~~~~ ~~~~~~~~~~~~~~ ~}\label{Chien:a1}\\
 {\varphi_{1}^{(C)}(0)=0,\;\;\;\;\varphi_{1}^{(C)}(1)=\frac{\lambda}{\lambda-1}\frac{d\varphi_{1}^{(C)}(y)}{dy}\bigg{|}_{y=1},~~~~~~~~}\\
  {-\int_{0}^{1}\frac{1}{\epsilon}\varphi_{1}^{(C)}(\epsilon)d\epsilon=\textbf{1};}\qquad \qquad \qquad \qquad \qquad~~~~~~~~
\end{split}
\right.
\end{equation}
\begin{equation}
 {W^{2}(0):} \left\{
\begin{split}
 {y^2\frac{d^{2}S_{1}^{(C)}(y)}{dy^{2}}=-\frac{1}{2}\left(\varphi_{1}^{(C)}(y)\right)^{2}, \;\;\;\;} ~~~~~~~~~~~~~~~~~~~~~~~~\\
 {S_{1}^{(C)}(0)=0,\;\;\;\;S_{1}^{(C)}(1)=\frac{\mu}{\mu-1}\frac{dS_{1}^{(C)}(y)}{dy}\bigg{|}_{y=1};~~~~~~~~~}
\end{split}
\right.
\end{equation}
$~~~~~\cdots \cdots $
\begin{equation}
 {W^{2i+1}(0):} \left\{
\begin{split}
 {y^2\frac{d^{2}\varphi_{i+1}^{(C)}(y)}{dy^{2}}=\sum_{j=1}^{i}\varphi_{j}^{(C)}(y)S_{i-j+1}^{(C)}(y)+Q_{i+1}^{(C)}y^{2}, } \;\;\;\;~~~~~~~\\
 {\varphi_{i+1}^{(C)}(0)=0,\;\;\varphi_{i+1}^{(C)}(1)=\frac{\lambda}{\lambda-1}\frac{d\varphi_{i+1}^{(C)}(y)}{dy}\bigg{|}_{y=1},~~~~~~~~~~~~~}\\
 {-\int_{0}^{1}\frac{1}{\epsilon}\varphi_{i+1}^{(C)}(\epsilon)d\epsilon=0;~~~~~~~~~~~~~~~~~~~~~~~~~~~~~~~~~~~~~~~~}
\end{split}
\right.
\end{equation}
\begin{equation}
 {W^{2i+2}(0):} \left\{
\begin{split}
 {y^2\frac{d^{2}S_{i+1}^{(C)}(y)}{dy^{2}}=-\frac{1}{2}\sum_{j=1}^{i+1}\varphi_{j}^{(C)}(y)\varphi_{i-j+2}^{(C)}(y), } ~~~~~~~~~~~~~~\; \\
 {S_{i+1}^{(C)}(0)=0,\;\;\;\;S_{i+1}^{(C)}(1)=\frac{\mu}{\mu-1}\frac{dS_{i+1}^{(C)}(y)}{dy}\bigg{|}_{y=1};~~~~~~~~~~}
\end{split}
\right. \label{Chien:s3}
\end{equation}
 {In the case of clamped boundary ($\lambda=0$, $\mu=20/7$), we have the Chien's perturbation solutions}
\begin{eqnarray}
 {\varphi_{1}^{(C)}(y)} &=&  {-2y+2y^2,} \label{Chien:phi1}\\
 {s_{1}^{(C)}(y)} &=&  {\frac{41y}{42}-y^2+\frac{2y^3}{3}-\frac{y^4}{6},}\\
 {\varphi_{2}^{(C)}(y)} &=&  {\frac{233y}{1890}-\frac{2179y^2}{3780}+\frac{83y^3}{126}-\frac{5y^4}{18}+\frac{y^5}{12}-\frac{y^6}{90},}\\
 {s_{2}^{(C)}(y)} &=&  {-\frac{211y}{19845}+\frac{233y^2}{1890}-\frac{529y^3}{2268}+\frac{667y^4}{3240}-\frac{59y^5}{630} } \nonumber\\
&\;\;& {+\frac{13y^6}{540}-\frac{17y^7}{3780}+\frac{y^8}{2520},}\label{Chien:s2}\\
\cdots \cdots&\;& \nonumber
\end{eqnarray}
Unfortunately,  as pointed out by Volmir \cite{Volmir},  the deflection curve  given  by the Chien's perturbation method \cite{Qian} becomes concave at centre when the central deformation increases to a certain level, which is obviously in contradiction with physical phenomena.  

In summary,  for circular plate under uniform external pressure,  the  previous  perturbation  methods \cite{Vincent, Qian}  are  valid  only  for the small physical parameters, i.e.  corresponding to the weak nonlinearity.  

In $1965$ a modified iteration method was proposed by Yeh and Liu \cite{YehLiu}, which inherits the merits of iteration technique and Chien's perturbation method \cite{Qian}.   The  procedures of the modified iteration method  \cite{YehLiu, Zheng} are as follows:
\begin{eqnarray}
    y^{2}\frac{d^{2}\psi_{n}(y)}{dy^{2}}&=&-\frac{1}{2}\vartheta_{n}^{2}(y),\label{modified2:original:U}\\
    y^{2}\frac{d^{2}\vartheta_{n+1}(y)}{dy^{2}}&=&\vartheta_{n}(y)\psi_{n}(y)+Q_{n}y^{2}, \label{modified1:original:U}
\end{eqnarray}
subject to the boundary conditions
\begin{equation}
    \vartheta_{n+1}(0)=\psi_{n}(0)=0,  \label{modified1:boundary:U}
\end{equation}
\begin{equation}
    \vartheta_{n+1}(1)=\frac{\lambda}{\lambda-1}\cdot\frac{d\vartheta_{n+1}(y)}{dy}\bigg{|}_{y=1},\;\;\;\; \psi_{n}(1)=\frac{\mu}{\mu-1}\cdot\frac{d\psi_{n}(y)}{dy}\bigg{|}_{y=1},  \label{modified2:boundary:U}
\end{equation}
with the restriction condition
\begin{equation}
    W(0)=a=-\int_{0}^{1}\frac{1}{\varepsilon}\vartheta_{n+1}(\varepsilon)d\varepsilon \label{modified3:boundary:U}
\end{equation}
and the initial guess
\begin{equation}
    \vartheta_{1}(y)=\frac{-2a}{2\lambda+1}[(\lambda+1)y-y^{2}].  \label{modified:initial:U}
\end{equation}
However, Zhou \cite{ZhouYH} studied the relationship between Chien's perturbation solutions \cite{Qian}  and the modified iterative solutions \cite{YehLiu}, but found that they have the {\em same} convergent region.  The modified iterative method \cite{YehLiu}, therefore, is also only valid for a small central deflection, too.  Therefore,  iteration itself can {\em not}  enlarge the convergence radius of perturbation series, although greatly boosts the computational efficiency.

 Keller and Reiss \cite{Keller}  proposed the interpolation iterative method by introducing  an interpolation parameter to the iteration procedure.   Fortunately,  the interpolation iterative method  yields convergent solutions even for  loads as large as $Q=7000$.    However,  while all of the iterations can be obtained explicitly as polynomials, their degrees increase geometrically.   Keller and Reiss \cite{Keller}, therefore, computed the iterations approximately by means of finite differences and gave numerical results.   In $1988$, Zheng and Zhou \cite{Zheng2}  {\em proved} that the series solutions given by the  interpolation iterative method are convergent for {\em arbitrary} values of load if {\em proper}  interpolation iterative parameter is chosen \cite{Zheng2}.  

In this paper, the same problem is  solved  by means of the homotopy analysis method (HAM) \cite{liaoPhd, Liaobook, liaobook2, KV2012}, an analytic approximation technique  proposed by Liao \cite{liaoPhd} for highly nonlinear problems.   Unlike perturbation technique, the HAM is independent of any small/large physical parameters.  Besides, it provides great freedom to choose equation-type and solution expression of the high-order approximation equations.  Especially, the HAM also provides us a convenient way to guarantee the convergence of series solutions  by means of introducing the so-called ``convergence-control parameter" $c_{0}$.   {It should be emphasized that some mathematical theorems of convergence have been rigorously proved in the frame of the HAM \cite{liaobook2}.  For instance, it has been proved \cite{liaobook2} that the power series given by the HAM 
\begin{equation}
u(t)=\lim_{m\rightarrow+\infty}\sum_{n=0}^{m}\mu_{0}^{m,n}(c_{0})(-t)^n, \label{intro:power:series}
\end{equation}
where
\begin{equation}
 \mu_{0}^{m,n}(c_{0})=(-c_{0})^{n}\sum_{k=0}^{m-n}\binom{n-1+k}{k}(1+c_{0})^{k},
\end{equation}
converge to $1/(1+t)$ in the intervals:
\begin{equation}
-1<t<-\frac{2}{c_{0}}-1,\qquad when \;\;c_{0}<0,
\end{equation}
and
\begin{equation}
-\frac{2}{c_{0}}-1<t<-1, \qquad when\;\;c_{0}>0,
\end{equation}
respectively. So, the power series (\ref{intro:power:series}) converges to $1/(1+t)$ either in the interval $(-1,+\infty)$ if $c_{0}<0$ impends $0$, or in the interval $(-\infty,-1)$ if $c_{0}>0$ tends to $0$, respectively. In other words, the introduction of the convergence-control parameter $c_{0}$ allows the power series (\ref{intro:power:series}) to converge to $1/(1+t)$ in its entire definition domain. Note that the traditional power series (regarding $t$ as a small variable):
\begin{equation}
\frac{1}{1+t}\sim 1-t+t^{2}-t^{3}+t^{4}-\cdots
\end{equation}
only converges in the interval $(-1,1)$. Thus, the so-called ``convergence-control parameter" $c_{0}$  can indeed greatly enlarge the convergence interval of solution series. } As a powerful technique to solve highly nonlinear equations, the HAM has been successfully employed to solve various types of nonlinear problems over the past two decades \cite{Abbasbandy2006,Hayat2006, KV2008, Liang2010, Ghotbi2011, Nassar2011, Mastroberardino2011, Aureli2014, Duarte2015, Nagarajaiah2015}.   {Note that the F\"{o}ppl-Von K{\'a}rm{\'a}n's plate equations  was  solved  by  means of the HAM \cite{RA2012}}.  Especially, the HAM  can bring us something completely new/different: the steady-state resonant waves were first predicted by the HAM in theory \cite{xu2012JFM, Liu2014JFM, Liu2015JFM, Liao2016JFM} and then confirmed experimentally in a laboratory   \cite{Liu2015JFM}.  

In this paper,  Von K{\'a}rm{\'a}n's plate equations in the differential form with clamped boundary are solved at first.  By means of the normal HAM (without iteration), convergent results are obtained in the case of   {$w(0)/h=3.03$}, which is larger than the maximum convergent range ($w(0)/h<2.44$) of the perturbation method \cite{Chen}.  Further, an  iteration approach  is  proposed in the frame of the HAM to gain convergent solutions within a  rather large ratio of $w(0)/h > 20$, corresponding to a case of rather high nonlinearity.  Our HAM approximations agree well with those given by the interpolation iterative method \cite{Zheng}.   Furthermore,  analytic approximations for other three boundaries (moveable clamped, simple support and simple hinged support) are also presented in a similar way.   In addition,  we  prove  that  the previous perturbation methods (including Vincent's \cite{Vincent} and Chien's \cite{Qian}  perturbation methods)  and the modified iteration method \cite{YehLiu} are only special cases of the HAM.     

\section{Analytic approach based on the HAM}

Like  Zheng \cite{Zheng},  we express $\varphi(y)$ and $S(y)$  in power series
\begin{equation}
   \varphi(y)=\sum_{m=1}^{+\infty} a_{m}\cdot y^{m},\;\;\;\;S(y)=\sum_{m=1}^{+\infty}b_{m}\cdot y^{m}, \label{homotopy:series}
\end{equation}
where $a_{m}$ and $b_{m}$ are constant coefficients to be determined. This provides us the so-called ``solution expression'' of $\varphi(y)$ and $S(y)$ in the frame of the HAM.  Writing  
\begin{equation}
    W(0)=a, \label{variable}
\end{equation}
we have due to Eq.~(\ref{deflection}) an algebraic equation:
\begin{equation}
    \int_{0}^{1}\frac{1}{\varepsilon}\varphi(\varepsilon)d\varepsilon=-a. \label{introduce:boundary}
\end{equation}

 Let $\varphi_{0}(y)$ and $S_{0}(y)$  denote  the initial guesses of $\varphi(y)$ and $S(y)$, respectively, which satisfy the boundary conditions (\ref{ic1:original:U}), (\ref{ic2:original:U}) and (\ref{introduce:boundary}).  Moreover,   let $\mathcal{L}$ denote an auxiliary linear operator with property ${\cal L}[0]=0$,  $H_{1}(y)$ and $H_{2}(y)$ the auxiliary functions,  $c_{0}$ a non-zero auxiliary parameter, called the convergence-control parameter, and $q\in[0,1]$ the embedding parameter, respectively.   We construct  a  family of  differential equations in $q\in[0,1]$:
\begin{eqnarray}
    &\;\;&(1-q)\mathcal{L}\big{[}\Phi(y;q)-\varphi_{0}(y)\big{]} \nonumber\\
    &=&c_{0} q H_{1}(y) \bigg{[}y^{2}\frac{\partial^{2}\Phi(y;q)}{\partial y^{2}}-\Phi(y;q)\Xi(y;q)-\Theta(q)y^{2}\bigg{]}, \label{zeroth1}\\
   &\;\;& (1-q)\mathcal{L}\big{[}\Xi(y;q)-S_{0}(y)\big{]} \nonumber \\
   &=&c_{0} q H_{2}(y) \bigg{[}y^{2}\frac{\partial^{2}\Xi(y;q)}{\partial y^{2}}+\frac{1}{2}\Phi^{2}(y;q)\bigg{]}, \label{zeroth2}
\end{eqnarray}
subject to the boundary conditions
\begin{equation}
    \Phi(0;q)=\Xi(0;q)=0,   \label{ic3:original:U}
\end{equation}
\begin{equation}
    \Phi(1;q)=\frac{\lambda}{\lambda-1}\cdot\frac{\partial\Phi(y;q)}{\partial y}\bigg{|}_{y=1},\;\;\;\;\Xi(1;q)=\frac{\mu}{\mu-1}\cdot\frac{\partial\Xi(y;q)}{\partial y}\bigg{|}_{y=1},  \label{ic4:original:U}
\end{equation}
with the restriction condition
\begin{equation}
    \int_{0}^{1}\frac{1}{\varepsilon}\Phi(\varepsilon;q)d\varepsilon=-a.    \label{ic5:original:U}
\end{equation}
Note that  $\Phi(y;q), \Xi(y;q)$ and $\Theta(q)$ correspond to the unknown $\varphi(y), S(y)$ and $Q$, respectively, as mentioned below.      

Note that $Q$ is unknown for a given value of $W(0)=a$.    Expand $\Theta(q)$ in a power series
\begin{equation}
    \Theta(q) =  Q_{0}+\sum_{m=1}^{+\infty}Q_{m} \; q^{m},\label{another:Q}
\end{equation}
where 
\begin{equation}
   Q_{m} = {\cal D}_{m}[\Theta(q)], \label{homotopy:Q}
\end{equation}
is determined later,  in which
\begin{equation}
   {\cal D}_{m}[f]=\frac{1}{m!}\frac{\partial^{m}f}{\partial q^{m}}\bigg{|}_{q=0} \label{Dm}
\end{equation}
is called the $m$th-order homotopy-derivative of $f$.    

When $q=0$, due to the property $\mathcal{L}(0)=0$, Eqs.~(\ref{zeroth1})-(\ref{ic5:original:U}) have the solution
\begin{equation}
    \Phi(y;0)=\varphi_{0}(y),\;\;\;\;\;\;\Xi(y;0)=S_{0}(y). \label{initial:condition}
\end{equation}
 When $q=1$, Eqs.~(\ref{zeroth1})-(\ref{ic5:original:U}) are equivalent to the original equations (\ref{geq1:original:U})-(\ref{ic2:original:U}) and (\ref{introduce:boundary}),  provided
\begin{equation}
    \Phi(y;1)=\varphi(y),\;\;\;\;\Xi(y;1)=S(y),\;\;\;\;\Theta(1)=Q. \label{embedding:parameter}
\end{equation}
Therefore, as the embedding parameter $q$ increases from $0$ to $1$, $\Phi(y;q)$ varies (or deforms) continuously from the initial guess $\varphi_{0}(y)$ to $\varphi(y)$,  so do $\Xi(y;q)$ from the initial guess $S_{0}(y)$ to $S(y)$,  and $\Theta(q)$ from $Q_{0}$ to $Q$, respectively.    So, we call Eqs.~(\ref{zeroth1})-(\ref{ic5:original:U}) the zeroth-order deformation equation.    {Note that $Q_{0}$ is an unknown constant at present, which will be determined later}.  

Using (\ref{initial:condition}), we have the power series 
\begin{eqnarray}
    \Phi(y;q)=\varphi_{0}(y)+\sum_{m=1}^{+\infty}\varphi_{m}(y) \; q^{m}, \;\;\;
    \Xi(y;q)  = S_{0}(y)+\sum_{m=1}^{+\infty}S_{m}(y) \; q^{m},\label{another:s}
\end{eqnarray}
in which
\begin{equation}
   \varphi_{m}(y) = {\cal D}_{m}[\Phi(y;q)],\;\;\;\;\;\;S_{m}(y) = {\cal D}_{m}[\Xi(y;q)], \label{homotopy:S}
\end{equation}
where ${\cal D}_{m}$ is defined by (\ref{Dm}).   

It is well-known that convergence radius of a power series  is  finite in general.   Fortunately, in the frame of the HAM,   we have great freedom to choose the auxiliary linear operator  $\cal L$ and especially the so-called convergence-control parameter $c_0$.  Assume that all of them are properly chosen so that the power series  (\ref{another:Q}) and (\ref{another:s}) are convergent at $q=1$.   Then,  according to (\ref{embedding:parameter}), we have the so-called homotopy-series solution
\begin{equation}
    \varphi(y) =\sum_{m=0}^{+\infty}\varphi_{m}(y),\;\;\;\;
    S(y) = \sum_{m=0}^{+\infty}S_{m}(y),\;\;\;\;
   Q   =  \sum_{m=0}^{+\infty}Q_{m}.
\end{equation}

The governing equations and boundary conditions of $\varphi_{m}(y)$, $S_{m}(y)$ and  {$Q_{m-1}$} are obtained in the following way.  Substituting (\ref{another:Q}) and (\ref{another:s}) into  {the zeroth-order  deformation equations} (\ref{zeroth1})-(\ref{ic5:original:U}) and then equating the like-power of $q$, we have the so-called $m$th-order deformation equations 
\begin{eqnarray}
   &\;\;& \mathcal{L}[\varphi_{m}(y)-\chi_{m}\varphi_{m-1}(y)] =c_{0}H_{1}(y) \delta_{1,m-1}(y), \label{highorder1}\\
    &\;\;&\mathcal{L}[S_{m}(y)-\chi_{m}S_{m-1}(y)]  = c_{0}H_{2}(y)\delta_{2,m-1}(y)
   \label{highorder2},
\end{eqnarray}
subject to the boundary conditions
\begin{equation}
    \varphi_{m}(0)=S_{m}(0)=0,  \label{boundary1}
\end{equation}
\begin{equation}
    \varphi_{m}(1)=\frac{\lambda}{\lambda-1}\cdot\frac{d\varphi_{m}(y)}{d y}\bigg{|}_{y=1},\;\;\;\;S_{m}(1)=\frac{\mu}{\mu-1}\cdot\frac{d S_{m}(y)}{d y}\bigg{|}_{y=1},  \label{boundary2}
\end{equation}
with the restriction condition
\begin{equation}
   \int_{0}^{1}\frac{1}{\varepsilon}\cdot\varphi_{m}(\varepsilon)d\varepsilon=0,  \label{boundary3}
\end{equation}
 {where }
\begin{equation}
\chi_m =\left\{
\begin{array}{cc}
 {0} & \mbox{when $m\leq 1$}, \\
 {1} & \mbox{when $m > 1$,}
\end{array}
\right. \label{def:chi}
\end{equation}
and
\begin{eqnarray}
\delta_{1,m-1}(y) &=&  y^{2}\frac{d^{2}\varphi_{m-1}(y)}{d y^{2}}-\sum_{k=0}^{m-1}\varphi_{k}(y)S_{m-1-k}(y)-Q_{m-1} y^{2}, \\
\delta_{2,m-1}(y) &=& y^{2}\frac{d^{2}S_{m-1}(y)}{d y^{2}}+\frac{1}{2}\sum_{k=0}^{m-1}\varphi_{k}(y)\varphi_{m-1-k}(y). 
\end{eqnarray}

 According to the solution expression (\ref{homotopy:series}), we choose
\begin{equation}
    \varphi_{0}(y)=\frac{-2a}{2\lambda+1}[(\lambda+1)y-y^{2}],\;\;\;\;\;\;\;\;S_{0}(y)  =  (\mu+1)y-y^{2} \label{initial:s}
\end{equation}
as the initial guesses of $\varphi(y)$ and $S(y)$, and such an  auxiliary linear operator
\begin{equation}
    \mathcal{L}[u(y)]=\frac{d^{2}u(y)}{dy^{2}}. \label{linear:operator}
\end{equation}
Note that the initial guesses (\ref{initial:s}) satisfy all boundary conditions.  Similarly, according to the solution expression (\ref{homotopy:series}), the auxiliary functions $H_{1}(y)$ and $H_{2}(y)$ should be properly chosen so as to make sure that the right-hand sides of the high-order deformation equations (\ref{highorder1}) and (\ref{highorder2}) are in the forms
\begin{eqnarray}
    c_{0}H_{1}(y)\delta_{1,m-1}(y)&=&\sum_{k=0}d_{1,k}\;y^{k}, \label{coefficient:ergodicity1}\\
    c_{0}H_{2}(y)\delta_{2,m-1}(y)&=&\sum_{k=0}d_{2,k}\; y^{k}, \label{coefficient:ergodicity2}
\end{eqnarray}
where $d_{1,k}$ and $d_{2,k}$ are constants.  The auxiliary function $H_{1}(y)$ and $H_{2}(y)$, therefore, must be in the form
\begin{equation}
    H_{1}(y)=H_{2}(y)=\frac{1}{y^{2}}. \label{auxiliary:function}
\end{equation}

Then,  we have the general solutions of the high-order deformation equation (\ref{highorder1}) and (\ref{highorder2}):  
\begin{equation}
    \varphi_{m}(y)=\chi_{m}\varphi_{m-1}(y)+c_{0}\int_{0}^{y}\int_{0}^{\eta}\frac{\delta_{1,m-1}(\tau)}{\tau^{2}}d\tau d\eta +D_{1,m}y+D_{2,m},  \label{solve:phi1}
\end{equation}
\begin{equation}
   S_{m}(y)=\chi_{m}S_{m-1}(y)+c_{0}\int_{0}^{y}\int_{0}^{\eta}\frac{\delta_{2,m-1}(\tau)}{\tau^{2}}d\tau d\eta +D_{3,m}y+D_{4,m} \label{solve:s1},
\end{equation}
 {where $D_{1,m}$, $D_{2,m}$, $D_{3,m}$, $D_{4,m}$  are  determined by four linear boundary conditions (\ref{boundary1}) and (\ref{boundary2}), and the unknown $Q_{m-1}$ is determined by the restriction condition (\ref{boundary3})}.  In this way, $\varphi_{m}(y)$ , $S_{m}(y)$ and $Q_{m-1}$ of Eqs.~(\ref{highorder1}) and (\ref{highorder2}) can be obtained step by step, starting from $m=1$.   Then, we have the $M$th-order homotopy-approximation of $\varphi(y)$, $S(y)$ and $Q$:
\begin{equation}
   \tilde{\varphi}(y)=\sum_{m=0}^{M}\varphi_{m}(y),\;\;\;\;
   \tilde{S}(y)=\sum_{m=0}^{M}S_{m}(y), \;\;\;\;\tilde{Q} =\sum_{m=0}^{M}Q_{m}.\label{solve:s}
\end{equation}

Define the  squared residual error
\begin{equation}
   {\cal E}=\int_0^1 \left\{\bigg{(}{\cal N}_{1}\left[\tilde{\varphi}(y),\tilde{S}(y),y\right]\bigg{)}^{2}+\bigg{(}{\cal N}_{2}\left[\tilde{\varphi}(y), \tilde{S}(y),y\right]\bigg{)}^{2}\right\} dy, \label{discrete:residual}
\end{equation}
where the nonlinear operators defined by ${\cal N}_1$ and ${\cal N}_2$ are related to the original equations (\ref{geq1:original:U}) and (\ref{geq2:original:U}).    {Obviouly, the smaller the $\cal E$, the more accurate the HAM approximation.}  {According to Liao \cite{Liaobook, liaobook2}, the optimal value of the convergence-control parameter $c_{0}$ is determined by the minimum of $\cal E$.  It was proved \cite{Liaobook, liaobook2} in general cases that the homotopy-series converge to solutions of original equations as long as all squared residual errors tend to zero.  So,  it is enough to check the squared residual error (\ref{discrete:residual}) only. }

According to Liao \cite{Liaobook}, convergence of the homotopy-series solutions can be greatly accelerated by means of iteration technique, which uses the $M$th-order homotopy-approximation (\ref{solve:s1})
as the new initial guesses $\varphi_{0}(y)$ and $S_{0}(y)$ for the next iteration.    This provides us the $M$th-order iteration of the HAM.   In the iteration process of the HAM,  the right-hand side of  Eqs.~(\ref{highorder1}) and (\ref{highorder2}) are truncated to $y^N$, say,
\begin{equation}
  c_{0}H_{1}(y)\cdot\delta_{1,m}(y)\approx \sum_{k=0}^{N}E_{m,k}\cdot y^{k}, \;\;\;\;
  c_{0}H_{2}(y)\cdot\delta_{2,m}(y)\approx \sum_{k=0}^{N}F_{m,k}\cdot y^{k}, \label{Nterms2}
\end{equation}
where $E_{m,k}$ and $F_{m,k}$ are constants and $N$ is called the truncation order.   

\section{Result analysis}

 Without loss of generality, the Poisson's ratio $\nu$ is taken to be $0.3$ in all cases considered in this paper.

\subsection{The HAM-based approach without iteration}

Without loss of generality, let us consider the clamped boundary in the case of  $a = 5$, equivalent to $w(0)/h = 3.03$.   At the beginning, the so-called convergence-control parameter $c_0$ is unknown.  Its optimal value (i.e. -0.28 in this case) is determined by the minimum of the squared residual error defined by (\ref{discrete:residual}).     
According to  Table~\ref{results:NoniterA5}, the squared residual error decreases to $4.9\times10^{-12}$ by means of $c_{0}=-0.28$.  So, unlike perturbation method that is valid  only for  $w(0)/h < 2.44$  in the same case  \cite{Chen},  we gain convergent result by means of the HAM.  This is mainly because the so-called convergence-control parameter in the frame of the HAM provides us a convenient way to guarantee the convergence of solution series.   

\begin{table}[tb]
\tabcolsep 0pt
\caption{The  approximations of $Q$, the squared residual errors $Err$ and the used CPU time in the case of $a = 5$ with the clamped boundary  by  means of  the HAM (without iteration) using $c_{0} = -0.28$.}
\vspace*{-12pt}\label{results:NoniterA5}
\begin{center}
\def\temptablewidth{0.8\textwidth}
{\rule{\temptablewidth}{1pt}}
\begin{tabular*}{\temptablewidth}{@{\extracolsep{\fill}}cccccc}
$m$, order of approx. &  $\cal E$    &  ~~~~~$Q$~~~~~  &CPU time (seconds)                      \\
\hline
20 &  $6.4\times10^{-2}$  &  131.7 & 6                          \\
40 &  $5.8\times10^{-4}$  &  132.1 & 35                  \\
60 &  $9.3\times10^{-6}$ &  132.1 & 106                   \\
80 &  $2.0\times10^{-7}$ &  132.2& 243                         \\
100 & $5.1\times10^{-9}$   &  132.2& 465                       \\
120 & $1.5\times10^{-10}$  &  132.2 & 782                      \\
140 & $4.9\times10^{-12}$&  132.2  & 1205                   \\
\end{tabular*}
{\rule{\temptablewidth}{1pt}}
 \end{center}
 \end{table}
 
 \begin{table}[tb]
\tabcolsep 0pt
\caption{The  approximations of $Q$, the squared residual errors $\cal E$ and the used CPU time versus iteration times in the case of $a = 5$  with the clamped boundary, given  by the1st-order HAM iteration approach using $c_{0} = -0.55$  with the truncation order $N=100$.}
\vspace*{-12pt}\label{results:IterA5}
\begin{center}
\def\temptablewidth{0.8\textwidth}
{\rule{\temptablewidth}{1pt}}
\begin{tabular*}{\temptablewidth}{@{\extracolsep{\fill}}cccccc}
$m$, times of iteration&$\cal E$ &~~~~$Q$ ~~~~  &CPU time (seconds)     \\
\hline
10& $4.5\times10^{-3}$ &  132.2 &5       \\
20& $4.8\times10^{-11}$&  132.2 &12      \\
30& $1.3\times10^{-18}$&  132.2 &18      \\
40& $1.7\times10^{-24}$ &  132.2 &25     \\
\end{tabular*}
{\rule{\temptablewidth}{1pt}}
 \end{center}
 \end{table}

Given a value of $a$,  the optimal value of $c_{0}$ for the Von K{\'a}rm{\'a}n's plate equations with clamped boundary can be obtained in a similarly, which can be expressed  by such an empirical formula
  {\begin{equation}
   c_{0} =-\frac{50}{50+0.45\cdot a^\frac{7}{2}}\;\;\;\;\;(0\leq a\leq5).  \label{noniter:parameter:c0}
\end{equation}}
 {Note that, like perturbation approximations,   all of these homotopy-approximations of $\varphi(y)$ and $S(y)$ are expressed in polynomial for a given value of  $a$.}

In Appendix A, we prove that the  previous perturbation methods for an arbitrary  perturbation quantity (including the Vincent's  \cite{Vincent} and Chien's \cite{Qian} methods) are only the special cases of the HAM when $c_0=-1$.  This also explains why the HAM can give convergent results for larger $w(0)/h$:  the convergence-control parameter $c_0$ plays an important role in the guarantee of convergence of solution series.

\subsection{Convergence acceleration by means of iteration}

Iteration can be naturally introduced into the frame of the HAM to greatly accelerate the convergence, as illustrated by Liao \cite{Liaobook, liaobook2}.   Without loss of generality, let us first consider here the same case of $a = 5$ with the clamped boundary.   It is found that the squared residual error  arrives its minimum at $c_{0} \approx -0.55$ by means of the HAM-based 1st-order iteration approach with the truncation order $N=100$.   As shown in Table~\ref{results:IterA5},  the squared residual error  quickly decreases to $1.7\times10^{-24}$ in only $25$ seconds, which is about $100$ times faster than the HAM approach without iteration.  Thus,  from the viewpoint of computational efficiency, the iteration HAM approach is used below in the subsequent part of this paper. 

\begin{figure}
    \begin{center}
        \begin{tabular}{cc}
            \includegraphics[width=2.5in]{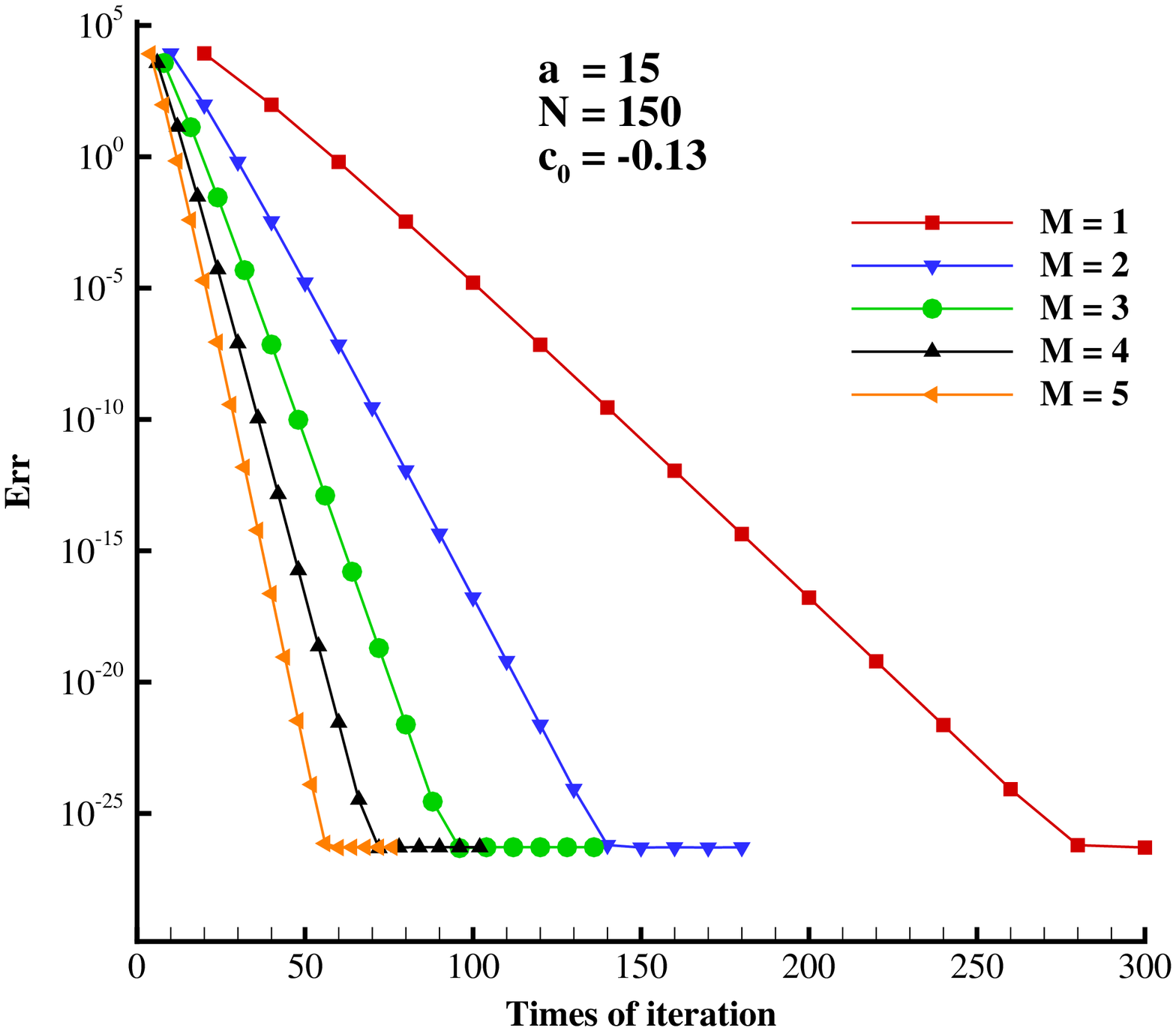} 
            \includegraphics[width=2.5in]{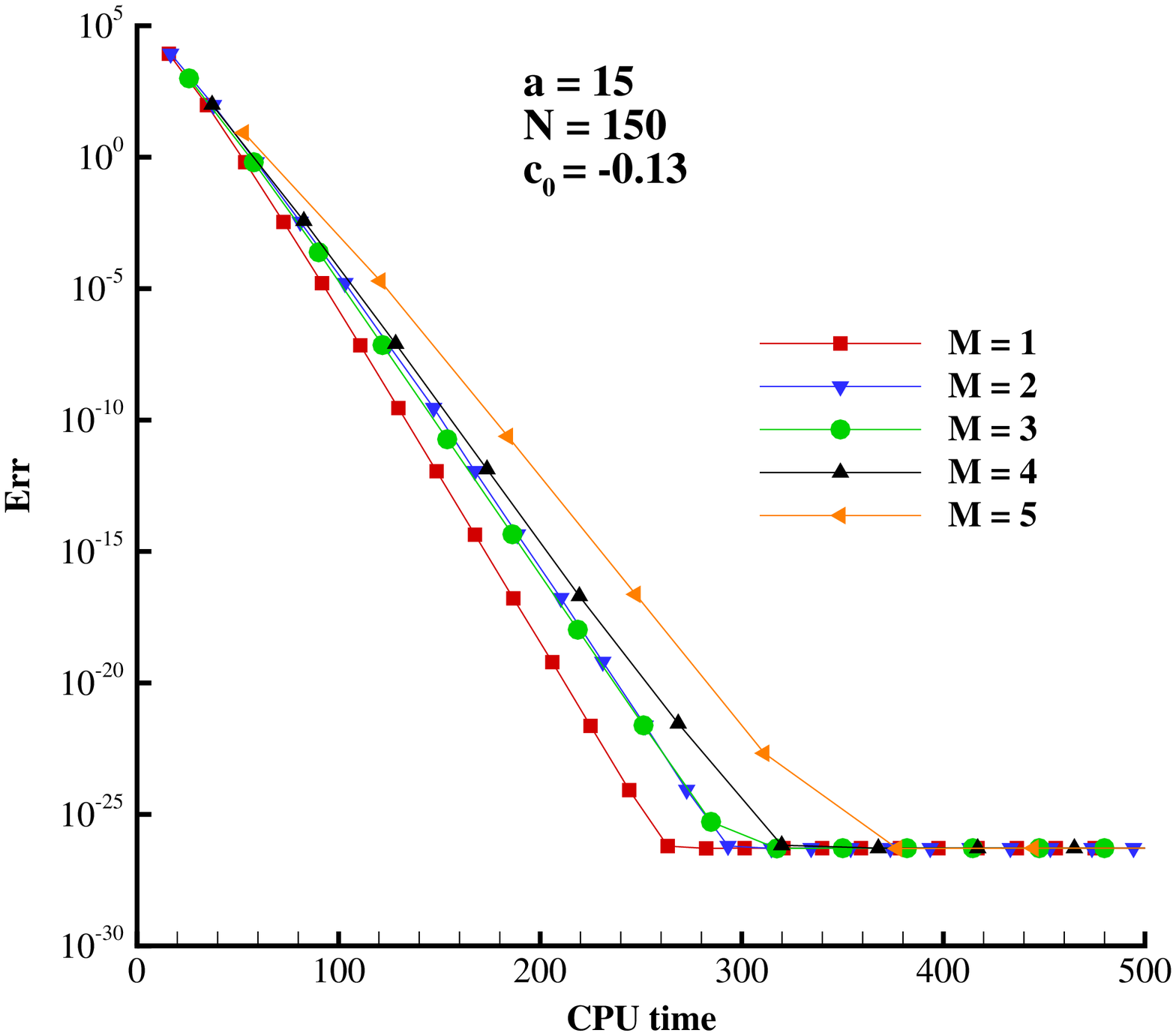} 
        \end{tabular}
    \caption{The squared residual error  versus the times of iteration and the CPU times in the case of $a=15$ with clamped boundary,  given by the HAM iteration approach using $c_{0}=-0.13$ with the truncation order $N=150$  and different order $M$ of iteration approach.  Square: 1st-order; Triangle down: 2nd-order; Circle: 3rd-order; Triangle up: 4th-order; Triangle left: 5th-order.} \label{figure:Iteration:times}
    \end{center}
\end{figure}

The HAM iteration approach contains two parameters, the iteration order $M$  and the truncation order $N$.  Without loss of generality, let us further consider the case of $a = 15$, corresponding to a higher nonlinearity.   As shown in Fig.~\ref{figure:Iteration:times}, the squared residual errors decrease to the level of $10^{-26}$ with different iteration order $M$:  the higher-order iteration need less iteration times, but more CPU times.  Thus, from the computational efficiency,  the first-order HAM iteration approach is suggested. 

 In addition, it is found that the iteration converges with high accuracy when $N$ is large enough, but larger $N$ also corresponds to a slower convergence.
Obviously,    as the load increases, i.e. $a$ enlarges, the nonlinearity becomes stronger, so that larger $N$ is necessary.   It is found that, in the frame of 1st-order HAM iteration approach, the  following empirical formula for the truncation order is good enough for all cases in this paper:  
\begin{equation}
 N = Max\left\{100, \gamma\cdot a\right\}, \label{def:gamma}
\end{equation}
where $\gamma$ is dependent upon the type of boundary:
\begin{equation}
   \gamma=\left\{
    \begin{array}{llllll}\nonumber
    10,~~~~~clamped~boundary;\\
    13,~~~~~moveable~clamped~boundary;\\
    ~7,~~~~~simple~support~boundary;\\
    ~5,~~~~~simple~hinged~support~boundary.
   \end{array}\right.
\end{equation}

\begin{table}[tb]
\tabcolsep 0pt
\caption{The homotopy-approximation of the load $Q$ versus $a$ for clamped boundary,  given by the HAM-based 1st-order iteration approach with the truncation order $N$ given by (\ref{def:gamma}) and the optimal convergence-control parameter $c_{0}$ given by (\ref{parameter:c0}).}
\vspace*{-12pt}\label{Clamped:Support}
\begin{center}
\def\temptablewidth{0.8\textwidth}
{\rule{\temptablewidth}{1pt}}
\begin{tabular*}{\temptablewidth}{@{\extracolsep{\fill}}cccccc}
$a$  &~~~~ $c_{0}$~~~~ &~~~~$N$ ~~~~ &  $Q$       \\
\hline
5    & -0.51 & 100  &  132.2        \\
10   & -0.21 &100   &  957.7      \\
15   & -0.10 & 150   &  3152.1     \\
20   & -0.06& 200   &   7386.9     \\
25   &-0.04 & 250   &   14334.1    \\
30   &-0.03 &300   &   24665.7    \\
35   &-0.02 & 350   &   39053.6      \\
\end{tabular*}
{\rule{\temptablewidth}{1pt}}
 \end{center}
 \end{table}
 
As shown in Table~\ref{Clamped:Support}, by means of the 1st-order HAM iteration approach,  the convergent  results in polynomials are obtained by means of the optimal convergence-control parameter with the empirical formula  
\begin{equation}
   c_{0} =-\frac{26}{26+a^2}\;\;\;\;\;(0\leq a\leq35) \label{parameter:c0}
\end{equation}
within $a\leq 35$, which is large enough for practice, since  $a = 35$ corresponds to $w(0)/h = 21.2$

\begin{figure}
    \begin{center}
        \begin{tabular}{cc}
            \includegraphics[width=2.5in]{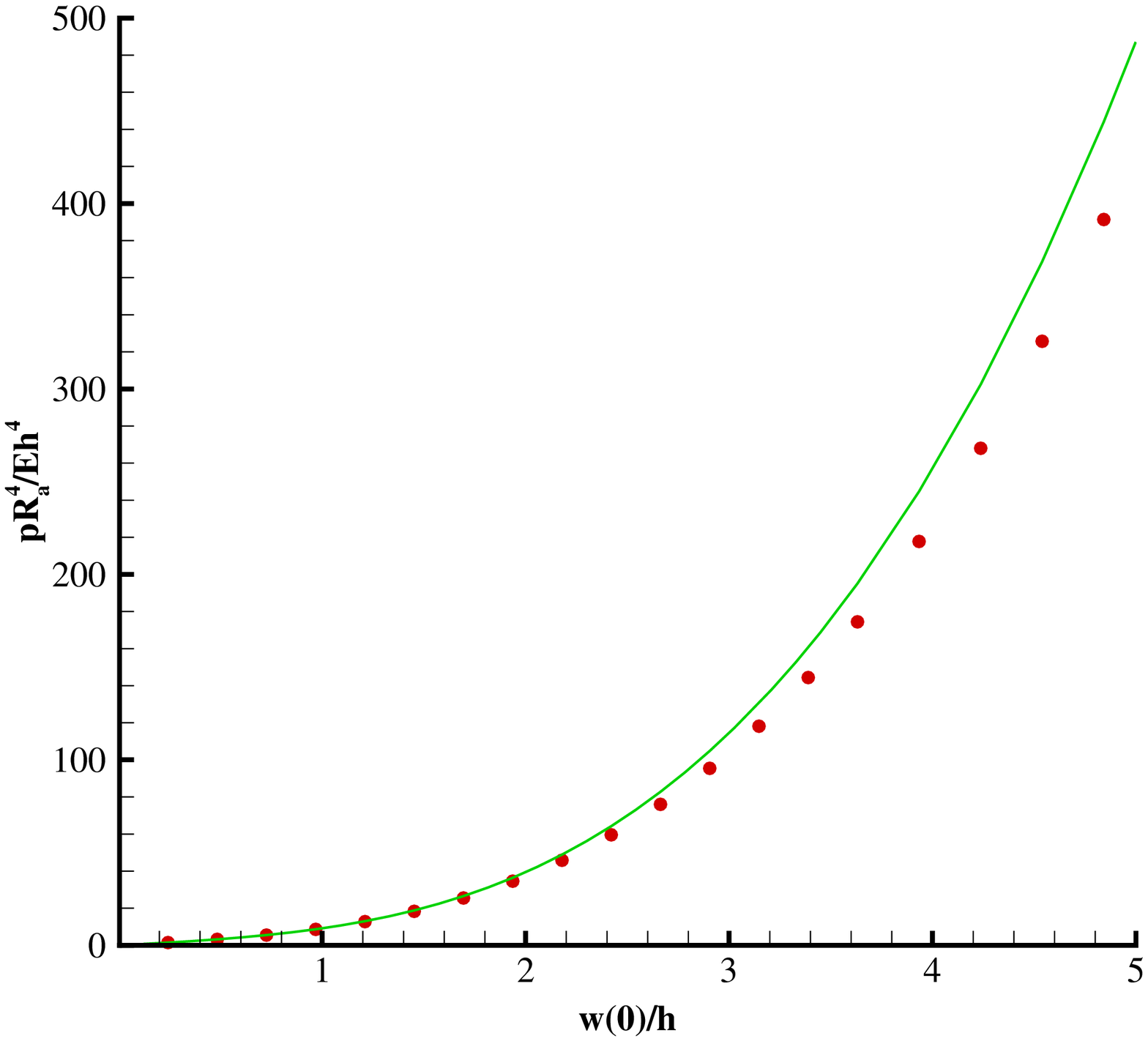} 
            \includegraphics[width=2.5in]{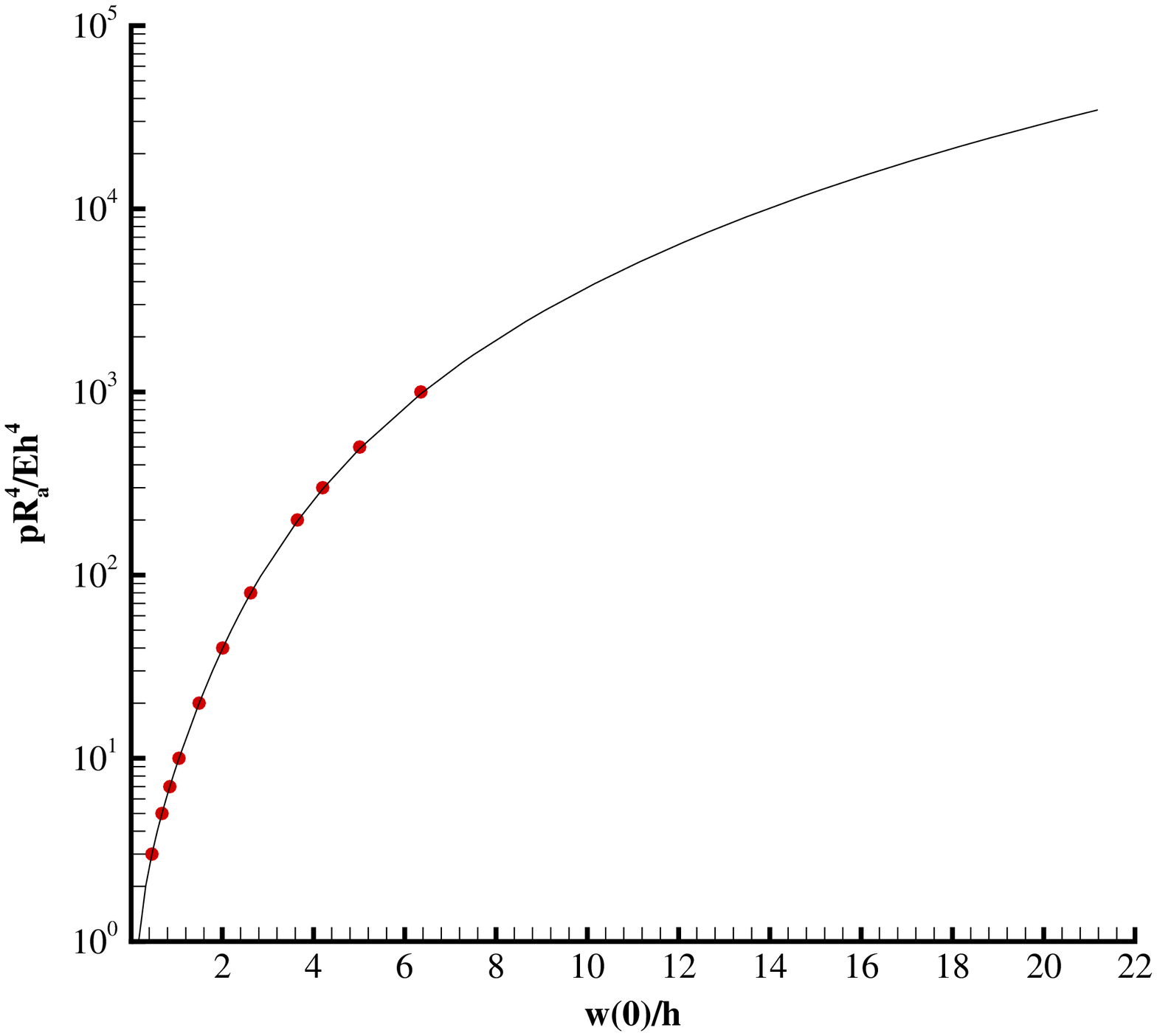} 
        \end{tabular}
    \caption{Comparison of the results given by the HAM-based 1st-order iteration approach and other methods for clamped boundary. Solid line: results given by the HAM; Symbols: results given by Chien's perturbation method \cite{Qian} (left) and by Zheng \cite{Zheng} using the interpolation iterative method  (right).} \label{Chien:Compare:HAM}
    \end{center}
\end{figure}

As shown in Fig.~\ref{Chien:Compare:HAM},  Chien's perturbation method \cite{Qian} is valid only in the region of $w(0)/h<2.44$, and becomes  worser and worser for larger $w(0)/h$.  In addition,  our HAM results agree quite well with those given by Zheng \cite{Zheng} using the interpolation iterative method, but converge in larger region.  The obtained deflections under different loads are depicted in Fig.~\ref{Clamped:Deform}.

In Appendix B,  we prove  that the modified iteration method \cite{YehLiu}  are only a special cases of the HAM when $c_0=-1$.  This reveals why the modified iteration method \cite{YehLiu} is valid for weak nonlinearity and besides shows the importance of the convergence-control parameter $c_0$ in the frame of the HAM.  In fact, it is the convergence-control parameter $c_0$ that differs the HAM from all other analytic approximations \cite{Liaobook, liaobook2}.   

\begin{figure}
    \begin{center}
        \begin{tabular}{cc}
            \includegraphics[width=2.5in]{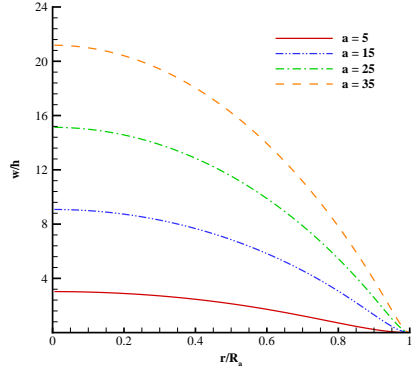} \\
        \end{tabular}
    \caption{The deflections of a thin circular plate with clamped boundary given by the HAM-based 1st-order iteration approach in the case of $a=5, 15, 25, 35$. Solid line: $pR_{a}^{4}/Eh^{4}=117.2$; Dash-double-dotted line: $pR_{a}^{4}/Eh^{4}=2795.2$; Dash-dotted line: $pR_{a}^{4}/Eh^{4}=12711.2$; Dashed line: $pR_{a}^{4}/Eh^{4}=34632.0$.} \label{Clamped:Deform}
    \end{center}
\end{figure}

Similarly, by means of the 1st-order HAM iteration approach,  the convergent  results for moveable clamped boundary are obtained by means of the optimal convergence-control parameter with the empirical formula
\begin{equation}
   c_{0} =-\frac{39}{39+a^2}  \label{c0:empirical:moveable}
\end{equation}
within the range of $a\leq 35$, as shown in Table~\ref{Moveable:Clamped}.   

For simple support boundary, the convergent  results are obtained by means of the optimal convergence-control parameter with the empirical formula
\begin{equation}
   c_{0} =-\frac{80}{80+a^2}  \label{c0:simple:support}
\end{equation}
within the range of $0\leq a\leq50$,  as shown in Table~\ref{Simple:Supports}.  For a circular plate with the boundary of simple hinged support,  the convergent  results are obtained by means of the optimal convergence-control parameter with the empirical formula
\begin{equation}
    c_{0} =-\frac{40}{40+a^\frac{5}{2}}  \label{c0:hinger:support}
\end{equation}
within the range of $0\leq a\leq50$, as shown in Table~\ref{Simple:Hinged}.

\begin{table}
\tabcolsep 0pt
\caption{The convergent  homotopy-approximation of the load $Q$ in case of different values of $a$ for a circular plate with moveable clamped boundary,  given by the HAM-based 1st-order iteration approach using the truncation order $N$ given by (\ref{def:gamma}) and the optimal convergence-control parameter $c_{0}$ given by (\ref{c0:empirical:moveable}).}
\vspace*{-12pt}\label{Moveable:Clamped}
\begin{center}
\def\temptablewidth{0.8\textwidth}
{\rule{\temptablewidth}{1pt}}
\begin{tabular*}{\temptablewidth}{@{\extracolsep{\fill}}cccccc}
$a$  &~~~~ $c_{0}$~~~~ &~~~~$N$ ~~~~ &  $Q$    \\
\hline
5    &-0.61 & 100 &   49.3  \\
10   &-0.28  &  130&   240.1  \\
15   & -0.15 & 195 &   657.7   \\
20   & -0.09 & 260 &   1372.5  \\
25   & -0.06 & 325 &   2450.9   \\
30   & -0.04 & 390 &  3956.8    \\
35   & -0.03 & 455 &   5952.2    \\
\end{tabular*}
{\rule{\temptablewidth}{1pt}}
 \end{center}
 \end{table}

\begin{table}
\tabcolsep 0pt
\caption{The convergent homotopy-approximation of the load $Q$ in case of different values of $a$ for a circular plate with the  simple support boundary,  given by the HAM-based 1st-order iteration approach with the truncation order $N$ given by (\ref{def:gamma}) and the optimal convergence-control parameter $c_{0}$ given by (\ref{c0:simple:support}).}
\vspace*{-12pt}\label{Simple:Supports}
\begin{center}
\def\temptablewidth{0.8\textwidth}
{\rule{\temptablewidth}{1pt}}
\begin{tabular*}{\temptablewidth}{@{\extracolsep{\fill}}cccccc}
$a$  &~~~~ $c_{0}$~~~~&~~~~$N$ ~~~~ &  $Q$  \\
\hline
10  &-0.44&100 &  107.8    \\
20  &-0.17&140 &  737.4    \\
30  &-0.08&210 &  2304.8    \\
40  &-0.05&280 &  5199.8  \\
50  &-0.03&350 &  9799.3    \\
\end{tabular*}
{\rule{\temptablewidth}{1pt}}
 \end{center}
 \end{table}

\begin{table}
\tabcolsep 0pt
\caption{The convergent homotopy-approximation of the load $Q$ in case of different values of $a$ for a circular plate with the boundary of simple hinged support given by the HAM-based 1st-order  iteration approach with the truncation order $N$ given by  (\ref{def:gamma}) and the optimal convergence-control parameter $c_{0}$ given by (\ref{c0:hinger:support}).}
\vspace*{-12pt}\label{Simple:Hinged}
\begin{center}
\def\temptablewidth{0.8\textwidth}
{\rule{\temptablewidth}{1pt}}
\begin{tabular*}{\temptablewidth}{@{\extracolsep{\fill}}cccccc}
$a$  &~~~~ $c_{0}$~~~~ &~~~~$N$ ~~~~ &  $Q$   \\
\hline
10   &-0.112&100&  890.0    \\
20  &-0.021& 100&  7152.3    \\
30   &-0.008&150&  24166.4    \\
40   &-0.004&200&  57308.7   \\
50  &-0.002&250 &  111955.3   \\
\end{tabular*}
{\rule{\temptablewidth}{1pt}}
 \end{center}
 \end{table}

\section{Conclusions}

In this paper, the homotopy analysis method (HAM) is applied to the large deflection of a circular thin plate under uniform external pressure.   By means of choosing a proper value of the so-called  convergence-control parameter $c_0$ given by the empirical formulas  (\ref{parameter:c0}), (\ref{c0:empirical:moveable}), (\ref{c0:simple:support}) and (\ref{c0:hinger:support}),   convergent results are successfully obtained even in the case of $w(0)/h > 20$ for four types of boundary conditions.   Besides, it is found that iteration can greatly accelerate the convergence of solutions.  In addition,  we reveal that the  previous perturbation methods for an arbitrary  perturbation quantity (including the Vincent's  \cite{Vincent} and Chien's \cite{Qian} methods) and the modified iteration method \cite{YehLiu}  are only the special cases of the HAM\footnote{Even the interpolation iterative method \cite{Keller} is also a special case of the HAM.  Limited to the length of the  paper, we will give the proof somewhere else} when $c_0=-1$.   This reveals the reason why the previous perturbation techniques \cite{Vincent,Qian}  and the modified iteration method \cite{YehLiu} are only valid for small physical parameters, corresponding to weak nonlinearity.   This work shows once again that the convergence-control parameter $c_0$ indeed plays a very important role in the frame of the HAM: it differs the HAM from all other analytic approximation methods.       

This paper demonstrates the validity of the HAM for the Von K{\'a}rm{\'a}n plate equations, and also clearly shows the superiority of the HAM over perturbation methods.   Without doubt, the HAM can be further applied to solve some challenging problems with high nonlinearity in solid mechanics.

\section*{Acknowledgment}
This work is partly supported by National Natural Science Foundation of China (Approval No. 11272209 and 11432009) and State Key Laboratory of Ocean Engineering (Approval No. GKZD010063).

\renewcommand{\theequation}{A-\arabic{equation}}
\setcounter{equation}{0}
\section*{Appendix A. Relations between the perturbation methods and the HAM approach}

Here  we  prove that  the perturbation methods for \emph{an arbitrary} perturbation quantity  (including Vincent's \cite{Vincent} and Chien's  \cite{Qian}  perturbation methods) are  special cases  of the HAM  approach when $c_0=-1$. 

 {In general, the perturbation solutions can be expressed as:}
\begin{equation}
 {u^{P}(y)=u^{P}_{0}(y)+u^{P}_{1}(y)\varepsilon+u^{P}_{2}(y)\varepsilon^{2}+\cdot\cdot\cdot, } \nonumber
\end{equation}
where $\varepsilon$ is a physical parameter.    
 {The HAM-series solutions are expressed by}
\begin{equation}
 {u^{H}(y;c_{0})=u^{H}_{0}(y;c_{0})+u^{H}_{1}(y;c_{0})+u^{H}_{2}(y;c_{0})+\cdot\cdot\cdot. } \nonumber
\end{equation}
 {Obviously,  if  $u^{P}_{i}(y)\varepsilon^{i}$ = $u^{H}_{i}(y;-1)$, then  $u^{P}(y)$ is exactly the same as  $u^{H}(y;-1)$,  say,   the perturbation method is a special case of the HAM when $c_{0}=-1$.}

First, we  describe the perturbation methods for a circular plate under uniform pressure in a general way.   Let $\zeta$ denote a perturbation quantity and write the perturbation series
\begin{equation}
\left\{ 
\begin{split}
    \varphi^{(P)}(y) =\sum_{i=1}^{+\infty}\varphi_{i}^{(P)}(y)\zeta^{2i-1}  ,\;\;\;
    S^{(P)}(y) = \sum_{i=1}^{+\infty}S_{i}^{(P)}(y)\zeta^{2i}, \;\;\;\;\\
    \;\;\;\;\;\;Q^{(P)} =\sum_{i=1}^{+\infty}Q_{i}^{(P)}\zeta^{2i-1}, \;\;\;\;\;
    W^{(P)}(0) =\sum_{i=1}^{+\infty}W_{i}^{(P)}(0)\zeta^{2i-1}.  \label{General:Expand:W}
\end{split} \right.
\end{equation}
Define
\begin{equation}
    S_{0}^{(P)}(y)=0. \nonumber
\end{equation}
Substituting (\ref{General:Expand:W}) into Eqs.~(\ref{geq1:original:U})-(\ref{ic2:original:U}) and (\ref{deflection}), and equating the like-power of $\zeta$, we have the governing equations
 \begin{eqnarray}
    y^{2}\frac{d^{2}\varphi_{m}^{(P)}(y)}{d y^{2}} & = &  \sum_{i=1}^{m}\varphi_{i}^{(P)}(y)S_{m-i}^{(P)}(y)+Q_{m}^{(P)}\cdot y^2,\label{General:highorder:phi}\\
    y^{2}\frac{d^{2}S_{m}^{(P)}(y)}{d y^{2}} & = & -\frac{1}{2}\sum_{i=1}^{m}\varphi_{i}^{(P)}(y)\varphi_{m+1-i}^{(P)}(y),
    \label{General:highorder:S}
\end{eqnarray}
subject to the boundary conditions
\begin{equation}
    \varphi_{m}^{(P)}(0)=S_{m}^{(P)}(0)=0,
    \label{General:boundary:1}
\end{equation}
\begin{equation}
    \varphi_{m}^{(P)}(1)=\frac{\lambda}{\lambda-1}\cdot\frac{d\varphi_{m}^{(P)}(y)}{dy}\bigg{|}_{y=1},\;\;\;\;\ S_{m}^{(P)}(1)=\frac{\mu}{\mu-1}\cdot\frac{dS_{m}^{(P)}(y)}{dy}\bigg{|}_{y=1},\label{General:boundary:2}
\end{equation}
with the restriction condition
\begin{equation}
    -\int_{0}^{1}\frac{1}{\varepsilon}\cdot \varphi_{m}^{(P)}(\varepsilon)d\varepsilon=W_{m}^{(P)}(0).
    \label{General:boundary:3}
\end{equation}
Note that there is an another \emph{linear} equation that characterizes the relation between the physical parameter $\zeta$ and the central deflection $W(0)$, which provides an another linear restriction condition. Note that the above-mentioned perturbation approach  is valid for {\em arbitrary} physical parameter $\zeta$.

The Von K{\'a}rm{\'a}n's equations for a circular plate under uniform pressure can be solved by the HAM in the following way, which is a little different from those mentioned in \S~2.

Let the initial guesses of $\varphi(y)$ and $S(y)$ be zero. We construct the following homotopy deformation equations:
\begin{eqnarray}
    (1-q)y^{2}\frac{d^{2}\tilde{\Phi}(y;q)}{dy^2}
    = c_{0} \bigg{[} q\; y^{2}\frac{d^{2} \tilde{\Phi}(y;q)}{dy^{2}}- \tilde{\Phi}(y;q)  \tilde{\Xi}(y;q)-\tilde{\Theta}(q)y^{2}\bigg{]}, \label{AppA:zeroth1}
\end{eqnarray}
\begin{eqnarray}
    (1-q)y^{2}\frac{d^{2}\tilde{\Xi}(y;q)}{dy^2}=c_{0}  \bigg{[}q\;y^{2}\frac{d^{2}\tilde{\Xi}(y;q)}{dy^{2}}+\left(\frac{1}{2q}\right)\tilde{\Phi}^{2}(y;q)\bigg{]}, ~~~~~~~~~~~ \label{AppA:zeroth2}
\end{eqnarray}
subject to the boundary conditions:
\begin{equation}
     \tilde{\Phi}(0;q)= \tilde{\Xi}(0;q)=0,   \label{AppA:left:boundary}
\end{equation}
\begin{equation}
     \tilde{\Phi}(1;q)=\frac{\lambda}{\lambda-1}\cdot\frac{\partial  \tilde{\Phi}(y;q)}{\partial y}\bigg{|}_{y=1},\;\; \tilde{\Xi}(1;q)=\frac{\mu}{\mu-1}\cdot\frac{\partial  \tilde{\Xi}(y;q)}{\partial y}\bigg{|}_{y=1},  \label{AppA:right:boundary}
\end{equation}
with the restriction condition:
\begin{equation}
    -\int_{0}^{1}\frac{1}{\varepsilon}\tilde{\Phi}(\varepsilon;q)d\varepsilon= \tilde{\Psi}(q),   \label{AppA:restrict:1}
\end{equation}
where $q\in[0,1]$ is the embedding parameter.    

When $q=0$,  the solutions of Eqs.~(\ref{AppA:zeroth1})-(\ref{AppA:restrict:1}) are the initial guess, i.e.
\begin{equation}
    \tilde{\Phi}(y;0) = 0,\;\;\;\;\;\;\;\; \tilde{\Xi}(y;0)=0.   \label{AppA:initial:solution}
\end{equation}
When $q=1$, Eqs.~(\ref{AppA:zeroth1})-(\ref{AppA:restrict:1}) are equivalent to the original equations (\ref{geq1:original:U})-(\ref{ic2:original:U}) and (\ref{introduce:boundary}), provided
\begin{equation}
    \tilde{\Phi}(y;1)=\varphi(y),\;\;\;\; \tilde{\Xi}(y;1)=S(y),\;\;\;\; \tilde{\Theta}(1)=Q,\;\;\;\; \tilde{\Psi}(1)=W(0). \label{AppA:provided}
\end{equation}
Then, according to (\ref{AppA:initial:solution}), $\tilde{\Phi}(y;q)$, $\tilde{\Xi}(y;q)$, $\tilde{\Theta}(q)$ and $\tilde{\Psi}(q)$ can be expanded as
\begin{equation}
\left\{ 
\begin{split}
    \tilde{\Phi}(y;q)=\sum_{m=1}^{+\infty}\tilde{\varphi}_{m}(y) \; q^{m},\;\;\;\;\;\;
    \tilde{\Xi}(y;q)  = \sum_{m=1}^{+\infty}\tilde{S}_{m}(y) \; q^{m},  \\
  \tilde{\Theta}(q) =\sum_{m=0}^{+\infty}\tilde{Q}_{m} \; q^{m},\;\;\;\;\;\;\;\;\;\;\;\;\;
    \tilde{\Psi}(q)  =  \sum_{m=0}^{+\infty}\tilde{W}_{m}(0) \; q^{m}.  \label{Appendix:B:series}
\end{split}
\right.
\end{equation}
 {Substituting  (\ref{Appendix:B:series}) into equations (\ref{AppA:zeroth1})-(\ref{AppA:restrict:1}) and equating the like-power of $q$, it is easy to find that $\tilde{Q}_{0}=0$ and $\tilde{W}_{0}(0)=0$,} and besides to obtain the $m$th-order deformation equations
\begin{eqnarray}
   &\;\;& y^2 \frac{d^2}{d y^2}\left[\tilde{\varphi}_{m}(y)-\chi_{m}\tilde{\varphi}_{m-1}(y)\right]\nonumber\\
    &=&c_{0} \; \bigg{(}y^{2}\frac{d^{2}\tilde{\varphi}_{m-1}(y)}{d y^{2}}-\sum_{i=1}^{m}\tilde{\varphi}_{i}(y)\tilde{S}_{m-i}(y)-\tilde{Q}_{m}y^{2}\bigg{)} , \label{AppA:highorder1}\\
    &\;\;& y^2 \frac{d^2}{d y^2}\left[\tilde{S}_{m}(y)-\chi_{m}\tilde{S}_{m-1}(y)\right]\nonumber\\
    &=&c_{0} \; \bigg{(}y^{2}\frac{d^{2}\tilde{S}_{m-1}(y)}{d y^{2}}+\frac{1}{2}\sum_{i=1}^{m}\tilde{\varphi}_{i}(y)\tilde{\varphi}_{m+1-i}(y)\bigg{)}.\label{AppA:highorder2}
\end{eqnarray}
subject to the boundary conditions
\begin{equation}
    \tilde{\varphi}_{m}(0)=\tilde{S}_{m}(0)=0,  \label{AppA:highorder:boundary1}
\end{equation}
\begin{equation}
    \tilde{\varphi}_{m}(1)=\frac{\lambda}{\lambda-1}\cdot\frac{d\tilde{\varphi}_{m}(y)}{d y}\bigg{|}_{y=1},\;\;\;\;\tilde{S}_{m}(1)=\frac{\mu}{\mu-1}\cdot\frac{d \tilde{S}_{m}(y)}{d y}\bigg{|}_{y=1},  \label{AppA:highorder:boundary2}
\end{equation}
with the restriction condition
\begin{equation}
   -\int_{0}^{1}\frac{1}{\varepsilon}\cdot \tilde{\varphi}_{m}(\varepsilon)d\varepsilon=\tilde{W}_{m}(0),  \label{AppA:highorder:boundary3}
\end{equation}
where $\chi_{m}$ is defined by (\ref{def:chi}).

 {When $c_0 = -1$,   Eqs.~(\ref{AppA:highorder1})-(\ref{AppA:highorder:boundary3}) are the same as the perturbation procedures (\ref{General:highorder:phi})-(\ref{General:boundary:3}), apart from the boundary conditions in the 1st-order deformation equations that leads to $\widetilde{\varphi}_{1}(y)=\zeta \cdot \varphi_{1}^{(P)}(y)$.  Thus, it holds}
\begin{equation}
\left\{ 
\begin{split}
 {\tilde{\varphi}_{m}(y) = \varphi_{m}^{(P)}(y)\; {\zeta^{2m-1}},   \;\;\;\;\;\; \tilde{S}_{m}(y)=S_{m}^{(P)}(y) \; {\zeta^{2m}},}\\
 {\tilde{Q}_{m} =Q_{m}^{(P)} \; {\zeta^{2m-1}},\;\;\;\;\;\;\;\;  \tilde{W}_{m}(0)=W_{m}^{(P)}(0) \; {\zeta^{2m-1}}.}
    \label{General:HAM:setQ}
\end{split} \right. 
\end{equation}
Therefore, the  perturbation  methods for {\em arbitrary}  perturbation quantity $\zeta$ are only special cases of the HAM when $c_{0}=-1$.  

 {For instance, if we choose the perturbation quantity $\zeta=Q$ and define   {$\widetilde{\Theta}(q)=Q\cdot q$ in the frame of the HAM, i.e. $\tilde{Q}_1=Q$ and $\tilde{Q}_m=0$ for $m\geq 2$},  we have}
 \begin{equation}
 {q:} \left\{
\begin{split}
 {y^2\frac{d^{2}\tilde{\varphi}_{1}(y)}{dy^{2}}=  {Q\cdot y^{2}}, \;\;\;\;y^2\frac{d^{2}\tilde{S}_{1}(y)}{dy^{2}}=-\frac{1}{2}\left(\tilde{\varphi}_{1}(y)\right)^{2},} ~~ \\
 {\tilde{\varphi}_{1}(0)=0,\;\;\;
\tilde{\varphi}_{1}(1)=\frac{\lambda}{\lambda-1}\frac{d\tilde{\varphi}_{1}(y)}{dy}\bigg{|}_{y=1},\;\;~~~~~~~~~~~~~}\\
 {\tilde{S}_{1}(0)=0,\;\;\;\;\tilde{S}_{1}(1)=\frac{\mu}{\mu-1}\frac{d\tilde{S}_{1}(y)}{dy}\bigg{|}_{y=1};~~~~~~~~~~~~~~}\label{AppA:q1}
\end{split}
\right.
\end{equation}

\begin{equation}
 {q^{2}:} \left\{
\begin{split}
 {y^2\frac{d^{2}\tilde{\varphi}_{2}(y)}{dy^{2}}=\tilde{\varphi}_{1}(y)\tilde{S}_{1}(y),~~~~~~~~~~~~~~~~~~~~~~~~~~~~~~~~~~ } \\
 {y^2\frac{d^{2}\tilde{S}_{2}(y)}{dy^{2}}=-\tilde{\varphi}_{1}(y)\tilde{\varphi}_{2}(y),~~~~~~~~~~~~~~~~~~~~~~~~~~~~~~~~ } \\
 {\tilde{\varphi}_{2}(0)=0,\;\;\tilde{\varphi}_{2}(1)=\frac{\lambda}{\lambda-1}\frac{d\tilde{\varphi}_{2}(y)}{dy}\bigg{|}_{y=1},\;\;~~~~~~~~~~~~~~~}\\
 {\tilde{S}_{2}(0)=0,\;\;\;\;\tilde{S}_{2}(1)=\frac{\mu}{\mu-1}\frac{d\tilde{S}_{2}(y)}{dy}\bigg{|}_{y=1};~~~~~~~~~~~~~~~}
\end{split}
\right.
\end{equation}
$~~~~~\cdots \cdots $
\begin{equation}
 {q^{i+1}:} \left\{
\begin{split}
 {y^2\frac{d^{2}\tilde{\varphi}_{i+1}(y)}{dy^{2}}=\sum_{j=1}^{i}\tilde{\varphi}_{j}(y)\tilde{S}_{i-j+1}(y), ~~~ ~~~~~~~~~~~~~~~~~~~~}\;\;\\
 {y^2\frac{d^{2}\tilde{S}_{i+1}(y)}{dy^{2}}=-\frac{1}{2}\sum_{j=1}^{i+1}\tilde{\varphi}_{j}(y)\tilde{\varphi}_{i-j+2}(y),~ ~~~~~~~~~~~~~~~~~~~ } \\
 {\tilde{\varphi}_{i+1}(0)=0,\;\;\tilde{\varphi}_{i+1}(1)=\frac{\lambda}{\lambda-1}\frac{d\tilde{\varphi}_{i+1}(y)}{dy}\bigg{|}_{y=1},\;\;~~~~~~~~~~}~\\
 {\tilde{S}_{i+1}(0)=0,\;\;\;\;\tilde{S}_{i+1}(1)=\frac{\mu}{\mu-1}\frac{d\tilde{S}_{i+1}(y)}{dy}\bigg{|}_{y=1};~~~~~~~~~~~~}\label{AppA:qi}
\end{split}
\right.
\end{equation}
 {In the case of circular plate with clamped boundary ($\lambda=0$, $\mu=20/7$), the homotopy solutions are:}
\begin{eqnarray}
 {\tilde{\varphi}_{1}(y)} &=&   {Q} {\cdot \left(-\frac{y}{2}+\frac{y^2}{2}\right), } \label{homotopy:vincent:phi1}\\
 {\tilde{S}_{1}(y)} &=&   {Q^{2}} {\cdot\left(\frac{41y}{672}-\frac{y^2}{16}+\frac{y^3}{24}-\frac{y^4}{96}\right),}\\
 {\tilde{\varphi}_{2}(y)} &=&   {Q^{3}} {\cdot\left(\frac{659y}{80640}-\frac{41y^2}{2688}+\frac{83y^3}{8064}-\frac{5y^4}{1152}+\frac{y^5}{768}-\frac{y^6}{5760}\right),}\\
 {\tilde{S}_{2}(y)} &=&   {Q^{4}} {\cdot\bigg{(}-\frac{2357y}{1505280}+\frac{659y^2}{322560}-\frac{1889y^3}{967680}+\frac{103y^4}{96768} ~~~~~~~~~~~~~~~ } \nonumber\\
&\;\;& {-\frac{59y^5}{161280}+\frac{13y^6}{138240}-\frac{17y^7}{967680}+\frac{y^8}{645120}\bigg{)},} \label{homotopy:vincent:S2}\\
\cdots \cdots&\;\;& \nonumber
\end{eqnarray}
 {which are exactly the same as (\ref{Vincent:phi1})-(\ref{Vincent:S2}) given by  the Vincent's perturbation method \cite{Vincent}.  Therefore,  the Vincent's perturbation method \cite{Vincent} is indeed a special case of the HAM when $c_0=-1$.}  

 {If we choose the perturbation quantity $\zeta=W(0)$, then $\tilde{\Psi}(q)=W(0)\cdot q$.   Similarly,  we have}
\begin{equation}
 {q:} \left\{
\begin{split}
 {y^2\frac{d^{2}\tilde{\varphi}_{1}(y)}{dy^{2}}=\widetilde{Q}_{1}y^{2}, \;\;\;\;y^2\frac{d^{2}\tilde{S}_{1}(y)}{dy^{2}}=-\frac{1}{2}\left(\tilde{\varphi}_{1}(y)\right)^{2},~~~~~~~~ } \\
 {\tilde{\varphi}_{1}(0)=0,\;\;
\tilde{\varphi}_{1}(1)=\frac{\lambda}{\lambda-1}\frac{d\tilde{\varphi}_{1}(y)}{dy}\bigg{|}_{y=1},\;\;~~~~~~~~~~~~~~~~~~~}\\
  {-\int_{0}^{1}\frac{1}{\epsilon}\tilde{\varphi}_{1}(\epsilon)d\epsilon=\textbf{W(0)}},~~~~~~~~~~~~~~~~~~~~~~~~~~~~~~~~~~~~~~\\
 {\tilde{S}_{1}(0)=0,\;\;\;\;\tilde{S}_{1}(1)=\frac{\mu}{\mu-1}\frac{d\tilde{S}_{1}(y)}{dy}\bigg{|}_{y=1};~~~~~~~~~~~~~~~~~~~~}
\end{split}
\right.
\end{equation}

\begin{equation}
 {q^{2}:} \left\{
\begin{split}
 {y^2\frac{d^{2}\tilde{\varphi}_{2}(y)}{dy^{2}}=\tilde{\varphi}_{1}(y)\tilde{S}_{1}(y)+\widetilde{Q}_{2}y^{2},~~~~~~~~~~~~~~~~~~~~~~~~~~~~ } \\
 {y^2\frac{d^{2}\tilde{S}_{2}(y)}{dy^{2}}=-\tilde{\varphi}_{1}(y)\tilde{\varphi}_{2}(y),~~~~~~~~~~~~~~~~~~~~~~~~~~~~~~~~~~~~ } \\
 {\tilde{\varphi}_{2}(0)=0,\;\;\tilde{\varphi}_{2}(1)=\frac{\lambda}{\lambda-1}\frac{d\tilde{\varphi}_{2}(y)}{dy}\bigg{|}_{y=1},\;\;~~~~~~~~~~~~~~~~~~~~}\\
 {-\int_{0}^{1}\frac{1}{\epsilon}\tilde{\varphi}_{2}(\epsilon)d\epsilon=0,~~~~~~~~~~~~~~~~~~~~~~~~~~~~~~~~~~~~~~~~~~~~}\\
 {\tilde{S}_{2}(0)=0,\;\;\;\;\tilde{S}_{2}(1)=\frac{\mu}{\mu-1}\frac{d\tilde{S}_{2}(y)}{dy}\bigg{|}_{y=1};~~~~~~~~~~~~~~~~~~~~}
\end{split}
\right.
\end{equation}
$~~~~~\cdots \cdots $
\begin{equation}
 {q^{i+1}:} \left\{
\begin{split}
 {y^2\frac{d^{2}\tilde{\varphi}_{i+1}(y)}{dy^{2}}=\sum_{j=1}^{i}\tilde{\varphi}_{j}(y)\tilde{S}_{i-j+1}(y)+\widetilde{Q}_{i+1}y^{2}, ~~~ ~~~~~~~~~~\;\;}\\
 {y^2\frac{d^{2}\tilde{S}_{i+1}(y)}{dy^{2}}=-\frac{1}{2}\sum_{j=1}^{i+1}\tilde{\varphi}_{j}(y)\tilde{\varphi}_{i-j+2}(y),~ ~~~~~~~~~~~~~~~~~~~~~ } \\
 {\tilde{\varphi}_{i+1}(0)=0,\;\;\tilde{\varphi}_{i+1}(1)=\frac{\lambda}{\lambda-1}\frac{d\tilde{\varphi}_{i+1}(y)}{dy}\bigg{|}_{y=1},\;\;~~~~~~~~~~~~~}\\
 {-\int_{0}^{1}\frac{1}{\epsilon}\tilde{\varphi}_{i+1}(\epsilon)d\epsilon=0,~~~~~~~~~~~~~~~~~~~~~~~~~~~~~~~~~~~~~~~~~~~~}\\
 {\tilde{S}_{i+1}(0)=0,\;\;\;\;\tilde{S}_{i+1}(1)=\frac{\mu}{\mu-1}\frac{d\tilde{S}_{i+1}(y)}{dy}\bigg{|}_{y=1};~~~~~~~~~~~~~~}
\end{split}
\right.
\end{equation}
 {In the case of circular plate with clamped boundary ($\lambda=0$, $\mu=20/7$), the homotopy solutions are:}
\begin{eqnarray}
 {\tilde{\varphi}_{1}(y)} &=&   {W(0)} {\cdot \left(-2y+2y^2 \right),} \label{homotopy:chien:phi1}\\
 {\tilde{S}_{1}(y)} &=&   {W^{2}(0)} { \cdot \left(\frac{41y}{42}-y^2+\frac{2y^3}{3}-\frac{y^4}{6} \right),} \\
 {\tilde{\varphi}_{2}(y)} &=&   {W^{3}(0)} { \cdot \left( \frac{233y}{1890}-\frac{2179y^2}{3780}+\frac{83y^3}{126}-\frac{5y^4}{18}+\frac{y^5}{12}-\frac{y^6}{90} \right),}\\
 {\tilde{S}_{2}(y)} &=&   {W^{4}(0)} { \cdot \bigg{(} -\frac{211y}{19845}+\frac{233y^2}{1890}-\frac{529y^3}{2268}+\frac{667y^4}{3240}-\frac{59y^5}{630} ~~~~~~~} \nonumber\\
&\;\;& {+\frac{13y^6}{540}-\frac{17y^7}{3780}+\frac{y^8}{2520}\bigg{)},} \\ \label{homotopy:Chien:s2}
\cdots \cdots &\;\;&  \nonumber
\end{eqnarray}
 {which are exactly the same as  (\ref{Chien:phi1})-(\ref{Chien:s2}) given by  Chien's perturbation method \cite{Qian}.   Therefore, the Chien's perturbation method \cite{Qian} is also a special case of the HAM when $c_0=-1$, too.}

However, it should be emphasized that the HAM provides us great freedom to choose the convergence-control parameter $c_0$.   As shown in \S~3,  by means of choosing proper values of $c_0$,  we gain convergent results in a range of $Q$ and $W(0)/h$ much larger than the known perturbation results.   This again illustrates the importance of the convergence-control parameter $c_0$ to the HAM.   In fact, it is the convergence-control parameter $c_0$ that differs the HAM from all of other analytic approximation techniques.   Mathematically, the above proof  also  reveals  the reason why the HAM has advantages over the perturbation methods.      

\renewcommand{\theequation}{B-\arabic{equation}}
\setcounter{equation}{0}

\section*{Appendix B.  Relations between the modified iteration method and the HAM-based iteration approach}

 {By means of the 1st-order HAM-based iteration approach, the new approximations $\varphi^*(y) = \varphi_0(y) + \varphi_1(y)$ and $S^*(y) = S_0(y) + S_1(y)$ are used as the new initial guesses $ \varphi_0(y), S_0(y) $  for  next iteration, since the HAM provides us the freedom to choose initial guesses, as illustrated for various types of nonlinear problems by Liao \cite{Liaobook,liaobook2}.}   According to (\ref{highorder1}) and (\ref{highorder2}), $\varphi_{1}(y)$ and $S_{1}(y)$ satisfy 
\begin{eqnarray}
    y^{2}\frac{d^{2}{\varphi}_{1}(y)}{d y^{2}} &=& c_0\left[y^{2}\frac{d^{2}{\varphi}_{0}(y)}{d y^{2}} -{\varphi}_{0}(y){S}_{0}(y)-{Q}_{0}y^{2}\right],  \label{Get:Phik:2}\\
    y^{2}\frac{d^{2}{S}_{1}(y)}{dy^{2}} &=&   c_0 \left[ y^{2}\frac{d^{2}{S}_{0}(y)}{dy^{2}}+\frac{1}{2}{\varphi}_{0}^{2}(y) \right].
\end{eqnarray}
Thus, $\varphi^{*}(y)$ and $S^{*}(y)$ are governed by 
\begin{eqnarray}
    y^{2}\frac{d^{2}{\varphi}^*(y)}{d y^{2}} &=&y^{2}\frac{d^{2}{\varphi}_{0}(y)}{d y^{2}} +y^{2}\frac{d^{2}{\varphi}_{1}(y)}{d y^{2}} \nonumber \\
     &=&(1+c_0)y^{2}\frac{d^{2}{\varphi}_{0}(y)}{d y^{2}} - c_0\left[{\varphi}_{0}(y){S}_{0}(y)+{Q}_{0}y^{2}\right],  \label{Get:Phik:2:A}\\
    y^{2}\frac{d^{2}{S}^*(y)}{dy^{2}} &=&  y^{2}\frac{d^{2}{S}_{0}(y)}{d y^{2}} +y^{2}\frac{d^{2}{S}_{1}(y)}{d y^{2}}\nonumber\\
     &=&(1+c_0)  y^{2}\frac{d^{2}{S}_{0}(y)}{dy^{2}}-\frac{c_0}{2}{\varphi}_{0}^{2}(y),\label{Get:Sk:2:A}
\end{eqnarray}
subject to the boundary conditions
 \begin{equation}
    {\varphi}^*(0)={S}^*(0)=0,\label{App:HAM:boundary1:A:2}
\end{equation}
\begin{equation}
    {\varphi}^*(1)=\frac{\lambda}{\lambda-1}\cdot\frac{d{\varphi}^*(y)}{dy}\bigg{|}_{y=1},\;\;\;\;\ {S}^*(1)=\frac{\mu}{\mu-1}\cdot\frac{d{S}^*(y)}{dy}\bigg{|}_{y=1},\label{App:HAM:boundary2:A:2}
\end{equation}
and the restriction condition
\begin{equation}
    W(0)=-\int_{0}^{1}\frac{1}{\varepsilon}{\varphi}^*(\varepsilon)d\varepsilon.\label{App:HAM:boundary3:A:2}
\end{equation}
In case of  $c_{0} = -1$,  we have
\begin{eqnarray}
    y^{2}\frac{d^{2}{\varphi}^*(y)}{d y^{2}} &=&  {\varphi}_{0}(y){S}_{0}(y)+{Q}_{0}y^{2},  \label{Get:Phik:2}\\
    y^{2}\frac{d^{2}{S}^*(y)}{dy^{2}} &=&  \frac{1}{2}{\varphi}_{0}^{2}(y),\label{Get:Sk:2}
\end{eqnarray}
subject to the boundary conditions
 \begin{equation}
    {\varphi}^*(0)={S}^*(0)=0,\label{App:HAM:boundary1}
\end{equation}
\begin{equation}
    {\varphi}^*(1)=\frac{\lambda}{\lambda-1}\cdot\frac{d{\varphi}^*(y)}{dy}\bigg{|}_{y=1},\;\;\;\;\ {S}^*(1)=\frac{\mu}{\mu-1}\cdot\frac{d{S}^*(y)}{dy}\bigg{|}_{y=1},\label{App:HAM:boundary2}
\end{equation}
and the restriction condition
\begin{equation}
    W(0)=-\int_{0}^{1}\frac{1}{\varepsilon}{\varphi}^*(\varepsilon)d\varepsilon.\label{App:HAM:boundary3}
\end{equation}
Since the initial guesses are given at the beginning,  we take the following iterative procedures:
\begin{enumerate}
\item[(1)]  Calculate ${S}^*(y)$ according to Eqs.~(\ref{Get:Sk:2})-(\ref{App:HAM:boundary2});
\item[(2)]  Replace ${S}_{0}(y)$  by ${S}^*(y)$ as the new initial guess, i.e. ${S}_{0}(y) = S^*(y)$;
\item[(3)]  Calculate $\varphi^*(y)$  and $Q_{0}$  according to Eqs.~(\ref{Get:Phik:2}) and (\ref{App:HAM:boundary1})-(\ref{App:HAM:boundary3});
\item[(4)]   Replace $\varphi_0(y)$  by  $\varphi^*(y)$  as the new initial guess, i.e. $\varphi_0(y) =\varphi^*(y)$.
\end{enumerate}
  In the $n$th times of iteration,  write 
\begin{equation}
  \Theta_{n}(y)={\varphi}^*(y),\;\;\;\;\;\Upsilon_{n-1}(y)={S}^*(y),\;\;\;\;\;F_{n-1}={Q}_0. \nonumber
\end{equation}
Then the HAM-based 1st-order iteration approach in case of $c_0=-1$  mentioned-above is expressed by
\begin{eqnarray}
    y^{2}\frac{d^{2}\Upsilon_{n-1}(y)}{dy^{2}} &=& -\frac{1}{2}\Theta_{n-1}^{2}(y),\label{App:HAM:origin1}\\
    y^{2}\frac{d^{2}\Theta_{n}(y)}{d y^{2}} &=& \Theta_{n-1}(y)\Upsilon_{n-1}(y)+F_{n-1}y^{2},\label{App:HAM:origin2}
\end{eqnarray}
subject to the boundary conditions
 \begin{equation}
    \Theta_{n}(0)=\Upsilon_{n-1}(0)=0,\label{App:HAM:boundarynew1}
\end{equation}
\begin{equation}
    \Theta_{n}(1)=\frac{\lambda}{\lambda-1}\cdot\frac{d\Theta_{n}(y)}{dy}\bigg{|}_{y=1},\;\;\;\;\Upsilon_{n-1}(1)=\frac{\mu}{\mu-1}\cdot\frac{d\Upsilon_{n-1}(y)}{dy}\bigg{|}_{y=1},\label{App:HAM:boundarynew2}
\end{equation}
and the restriction condition
\begin{equation}
    W(0)=a=-\int_{0}^{1}\frac{1}{\varepsilon}\Theta_{n}(\varepsilon)d\varepsilon.\label{App:HAM:boundarynew3}
\end{equation}
If we choose the initial guess 
\begin{equation}
    \Theta_{0}(y)=\frac{-2a}{2\lambda+1}[(\lambda+1)y-y^{2}],  \label{App:HAM:initial}
\end{equation}
then Eqs.~(\ref{App:HAM:origin1})-(\ref{App:HAM:initial}) are exactly the same as Eqs.~(\ref{modified2:original:U})-(\ref{modified:initial:U}) for the modified iteration method \cite{YehLiu, Zheng}.  Therefore,  the modified iteration method \cite{YehLiu, Zheng} is indeed a special case of the HAM-based 1st-order iteration approach when $c_{0}=-1$.

\section*{References}
\bibliographystyle{elsarticle-num}
\bibliography{plate}
\end{document}